\documentclass[11pt]{article}

\usepackage{amsfonts,amssymb,amsbsy, amsmath, dsfont}
\usepackage[german,english]{babel}
\usepackage[utf8]{inputenc}
\usepackage{graphicx}
\usepackage{epsfig}
\usepackage{gauss}
\usepackage{color}
\usepackage{amsthm}
\usepackage{framed}
\usepackage{enumerate}
\usepackage{tikz}
\usepackage{mathtools}

\usepackage[hyperindex,backref]{hyperref} 

\newcommand{\E}{{\mathbb E}}

\newcommand{\N}{{\mathbb N}}

\renewcommand{\P}{{\mathbb P}}
\newcommand{\Q}{{\mathbb Q}}
\newcommand{\R}{{\mathbb R}}

\newcommand{\Z}{{\mathbb Z}}

\newcommand{\cA}{{\cal A}}
\newcommand{\cB}{{\cal B}}
\newcommand{\cC}{{\cal C}}

\newcommand{\cG}{{\cal G}}

\newcommand{\cI}{{\cal I}}

\newcommand{\cL}{{\cal L}}

\newcommand{\cN}{{\cal N}}

\newcommand{\cV}{{\cal V}}

\newcommand{\cZ}{{\cal Z}}

\newcommand{\ind}{\mathds{1}}

\newcommand{\tofor}[1]{\overset{#1 \to \infty}{\longrightarrow}}

\DeclareMathOperator{\Ber}{Ber}
\DeclareMathOperator{\diam}{diam}

\numberwithin{equation}{section}	

\newtheoremstyle{thm}
	{8pt}
	{8pt}
	{\itshape}
	{}
	{\bfseries }
	{}
	{\newline}
	{}%

\newtheoremstyle{namedthm}
	{8pt}
	{8pt}
	{\itshape}
	{}
	{\bfseries }
	{}
	{5pt}
	{}%
	
\newtheoremstyle{def}
	{8pt}
	{8pt}
	{}
	{}
	{\bfseries }
	{}
	{\newline}
	{}%

\newtheoremstyle{nameddef}
	{8pt}
	{8pt}
	{}
	{}
	{\bfseries }
	{}
	{5pt}
	{}%
	
\swapnumbers
\theoremstyle{thm}
\newtheorem{Thm}{Theorem}[section]
\newtheorem{Lemma}[Thm]{Lemma}

\newtheorem{Cor}[Thm]{Corollary}

\theoremstyle{namedthm}

\newtheorem{NamedProp}[Thm]{Proposition}

\theoremstyle{def}
\newtheorem{Def}[Thm]{Definition}

\theoremstyle{nameddef}

\newlength{\XWidth}

\makeatletter
\def\moverlay{\mathpalette\mov@rlay}
\def\mov@rlay#1#2{\leavevmode\vtop{%
   \baselineskip\z@skip \lineskiplimit-\maxdimen
   \ialign{\hfil$\m@th#1##$\hfil\cr#2\crcr}}}
\newcommand{\charfusion}[3][\mathord]{
    #1{\ifx#1\mathop\vphantom{#2}\fi
        \mathpalette\mov@rlay{#2\cr#3}
      }
    \ifx#1\mathop\expandafter\displaylimits\fi}
\makeatother

\begin{document}
\author{Sebastian Ziesche\thanks{Karlsruhe Institute of Technology, sebastian.ziesche@kit.edu}
}

\title{Bernoulli Percolation on random Tessellations}
\date{\today}
\maketitle

\begin{abstract}
\noindent We generalize the standard site percolation model on the $d$-dimensional lattice to a model on random tessellations of $\R^d$. We prove the uniqueness of the infinite cluster by adapting the Burton-Keane argument \cite{burton1989density}, develop two frameworks that imply the non-triviality of the phase transition and show that large classes of random tessellations fit into one of these frameworks. Our focus is on a very general approach that goes well beyond the typical Poisson driven models. The most interesting examples might be Voronoi tessellations induced by determinantal processes or certain classes of Gibbs processes introduced in \cite{schreiber2013}. In a second paper we will investigate first passage percolation on random tessellations.
\end{abstract}

\maketitle

\begin{flushleft}
\end{flushleft}

\section{Introduction}\label{secintro}
The percolation model was introduced by Broadbend and Hammersley on $d$-dimensional lattices in the late fifties. In the meantime it was generalized to transitive or even quasi-transitive graphs and also to the random graph induced by the Poisson-Voronoi tessellation \cite{bollobas2006critical}. While the result of Bollobàs and Riordan, that the critical value for face percolation on the $2$-dimensional Poisson-Voronoi tessellation is excellent, it relies very much on the specifics of the Poisson-process and on the geometry of the plane. We want to start a very general investigation of Bernoulli face percolation on random tessellations of $\R^d$. By that, we mean a two-step model, where a tessellation is constructed randomly in the first step while in the second step each cell of this tessellation is colored black independently of each other cell with probability $p \in [0,1]$.\\
The structure of the paper is as follows. Section \ref{sec:perc_on_non-trans_graphs} contains generalizations of three well known theorems on percolation from transitive graphs to graphs with significantly less structure or symmetry. We show first, that $p_c > 0$, second that $p_c \geq \tfrac{1}{2}$ in the planar case (using the argument of Zhang) and third the uniqueness of the infinite cluster (adapting the proof of Gandolfi, Grimmett and Russo from \cite{gandolfi1988uniqueness2}). In these theorems the transitivity has been substituted by weaker properties that we will almost surely find in the realisations of our random tessellations. However these theorems are interesting in their own right.

In section \ref{sec:perc_on_random_tess} we proof the uniqueness of the infinite cluster by adapting the Burton Keane argument \cite{burton1989density} to our model. This will be possible under the extremely weak assumptions of stationarity and a moment condition on the cell-distribution. As a corollary we will show, that $p_c \geq \tfrac{1}{2}$ in the planar case. At this point one might hope, that an ergodic or maybe a mixing tessellation also exhibits a non-trivial phase transition. We will sketch a mixing counterexample that shows, that this is not the case in general.

Following this counterexample we will propose two frameworks in section \ref{sec:scale-mixing_tess} and \ref{sec:tame_tess}, that imply a non-trivial phase-transition. The first one is a kind of mixing condition. The second one is defined by auxiliary random fields, that encode how often very small or very large cells are observed. In both sections we will show for various classes of point processes, that the Voronoi tessellation induced by them fit into one of the two frameworks. This will include determinantal processes and certain classes of Poisson cluster and Gibbsian point processes introduced in \cite{schreiber2013}.

The emphasis in this paper lies on the large generality and basic results.

\section{Percolation on non-transitive graphs}\label{sec:perc_on_non-trans_graphs}
The purpose of this section is, to generalize the proofs of well known results from percolation on lattices. We want to do this in such a way, that the assumptions hold a.s. for various types random tessellations. We introduce the necessary notation first.

We work on an abstract probability space $(\Omega, \mathfrak{F}, \P)$. Let $I$ be a finite or countably infinite index set and $X := (X_i)_{i \in I}$ a family of i.i.d. random variables, where $X_i \sim \text{Ber}(p)$, $p \in [0,1]$ for all $i \in I$. We write $\P_p$ for the distribution of $X$. The space $\{0,1\}^I$ with the usual product $\sigma$-algebra is equipped with the canonical partial order $\preceq$, i.e.\ for any $\omega, \omega' \in \{0,1\}^I$ we have $\omega \preceq \omega'$ iff $\omega_i \leq \omega_i'$ for all $i \in I$. We call a real-valued function $f:\{0,1\}^I \to \R$ increasing, iff $f(\omega) \leq f(\omega')$ for all $\omega \preceq \omega'$. A function $f$ is decreasing iff $-f$ is increasing. An set $A \subset \{0,1\}^I$ is increasing (decreasing) iff the indicator function $\ind_A$ is increasing (decreasing). The well known FKG-inequality states, that for two increasing functions $f, g$
\[
  \E_p[f(X) g(X) ] \geq \E_{p}[f(X)] \E_{p}[g(X)],
\]
where $\E_{p}$ is the expectation with respect to $\P_p$. A proof can be found in the book of Grimmett \cite{grimmett1999percolation} which might as well be the best starting point for an introduction to percolation.

Without further mentioning, we always work on undirected, locally finite graphs without loops or double edges. Let $\cG = (V, E)$ be a graph with vertex set $V$ and edge set $E$. We say a set $U \subset V$ induces the subgraph $\cG'$ iff $\cG' = (U, \{\{v,w\} \in E \mid v, w \in U\})$. The site percolation model studies the family of Ber$(p)$ distributed random variables $X = (X_v)_{v \in V}$ and the random subgraph of $\cG$ induced by the set $\{v \in V \mid X_v = 1\}$.

A path $\gamma := (\gamma_1, \gamma_2, \dots)$, $\gamma_i \in V$ of $\cG$ is a finite or countably infinite sequence of vertices such that $\{\gamma_i, \gamma_{i+1}\} \in E$ for any $i \in \N$ and $\gamma_i \neq \gamma_j$ if $i \neq j$. We call the connected components of $\cG$ the clusters. In a cluster any two vertices are connected by a path in $\cG$. A cycle in $\cG$ is a finite path where starting and ending vertex are adjacent. We use $d_\cG$ for the usual graph metric and extend this notion in the way that $d_\cG(U, W) := \inf\{d_{\cG}(u, w) \mid u \in U, w \in W\}$ for $U, W \subset V$. For $v \in V, n \in \N$ the ball $B_n(v)$ is defined as $\{w \in V \mid d_{\cG}(w,v) \leq n\}$, and we write $B_n$ instead of $B_n(\mathbf{0})$, where $\mathbf{0}$ is the root of $\cG$. The outer and inner vertex boundary of a set $W \subset V$ is defined by $\partial^+ W := \{v \in V\setminus W \mid d_{\cG}(v,W) = 1\}$ and $\partial^- W := \{w \in W \mid d_{\cG}(w,V\setminus W) = 1\}$.

Percolation may be understood as the study of a randomly colored graph. To be precise, we call a pair $(\cG, c)$ a coloring of $\cG$, if $c$ is a map from $V$ to $\{0,1\}$. By a black path (cluster) we mean a path (cluster) of the subgraph of $\cG$ that is induced by the vertices $v$ with $c(v) = 1$. We write $C(v;\cG,c)$ for the black cluster that contains $v$ and remark that this might be the empty set if $v$ is white. In most cases we abbreviate this by $C_v$ if $\cG$ and $c$ are clear. In this section we will study the randomly colored graph $(\cG, X)$.

As usual
\[
  \theta_p(v) := \P_{p}[|C_v| = \infty ], \quad v \in \cG
\]
is the percolation function and
\begin{align}\label{eq:def:critical_prob}
  p_c := p_c(\cG) := \sup\{p \in [0,1] \mid \theta_p(\mathbf 0) = 0\}
\end{align}
is the critical value.

An obvious approach to proof results for percolation on random graphs is, to find conditions that hold for a.e.\ graph and study, what kind of results are implied by them. The most important property for many results in the classical theory on lattices is the transitivity of the underlying graph. This is completely destroyed in our model and we have to find substitutes for it.

If we take a look at the proof that $p_c > 0$ on a lattice, then we see, that the argument follows from an exponential growth condition on the number $b_{n}$ of paths of length $n \in \N$ starting in the origin. To be more precise, there is a constant $c \in \R$ depending on the lattice, such that $b_n \leq c^n$ for all natural number $n \in \N$. Apparently this condition won't hold a.s.\ for most random graphs stemming from random tessellations, as it is very often possible to observe cells with an arbitrary number of neighbouring cells. However the following trivial Lemma shows that this property can be weakened. This weakened property holds a.s.\ in framework II introduced in Section \ref{sec:tame_tess}. Let $\cA_n(\cG)$, $n \in \N$ be the set of connected subsets of size $n$ of $V$ that contain $\mathbf{0}$.

\begin{Lemma}\label{le:animal-bound-implies-p_c>0}
  Let $\cG$ be an infinite connected graph. If there is a $c \in \R$ such that $|\cA_n(\cG)| \leq c^n$ for all large enough $n \in \N$, then
  \[
    p_c(\cG) \geq \frac{1}{c}.
  \]
\end{Lemma}
\emph{Proof:} We have for all $p \in [0,1]$ and $n$ large enough
\[
  \theta_p(\mathbf 0) \leq \P_{p}\bigg[ \bigcup_{\alpha \in \cA_n(\cG)} \bigcap_{v \in \alpha} \{X_v = 1\}  \bigg] \leq \sum_{\alpha \in \cA_n(\cG)} p^n \leq (c p)^n.
\]
This implies $\theta_p(\mathbf{0}) = 0$ if $p < c^{-1}$.\qed
\bigskip

\subsection*{Uniqueness of the infinite cluster}
Another well known result is the uniqueness of the infinite cluster. This was first proven by Aizenman et. al. in \cite{aizenman1987uniqueness} and simplified by Gandolfi et. al. in \cite{gandolfi1988uniqueness2}. A short time afterwards Burton and Keane gave a very elegant new proof \cite{burton1989density}. While the first proof is a bit more technical and not as robust as the one of Burton and Keane, it doesn't rely as much on the transitivity and it can be quantified \cite{cerf2015lowerbound}. We will now show how the transitivity in \cite{gandolfi1988uniqueness2} can be relaxed, though we have to remark, that this generalization is maybe more of theoretical interest as the Burton Keane argument can be extended directly to percolation on random tessellations under extremely weak assumptions.

For any infinite connected graph $\cG = (V, E)$ and $v \in V$ let $L_v$ be the event, that $v$ is adjacent to two infinite black clusters in $(\cG,X)$.

\begin{Thm}\label{thm:uniqueness_GGR}
  Let $\cG = (V, E)$ be an infinite connected graph with the following properties:
  \begin{enumerate}
    \item The limit
        \[
          c_i := \lim_{n\to \infty } \frac{1}{|B_n|} \sum_{v \in B_n} g_i(v), \quad i \in \{1,2\}
        \]
        exists for the functions
        \[
          g_1:V \to \R:\ v \mapsto \P_p[L_v]
        \]
        and
        \[
          g_2:V \to \R:\ v \mapsto |\partial^+\{v\}|.
        \]
        If $c_1 = 0$ then $g_1(v) = 0$ for all $v \in V$.
    \item There is a $c_3 \in \R$ such that
        \[
          \lim_{n\to \infty } \frac{|B_n|}{n^{c_3}} = 0 \quad\text{and}\quad  \lim_{n\to \infty } \frac{|B_{n-\sqrt{n}}|}{|B_n|} = 1.
        \]
  \end{enumerate}
  Then there is at most one infinite cluster in $\cG$ for any $p \in [0,1]$.
\end{Thm}
\emph{Proof:} The proof is only a minor generalization of the one in \cite{gandolfi1988uniqueness2} but we do it for the sake of completeness. The statement is trivial for $p \in \{0,1\}$, hence let $p \in (0,1)$. If there are two infinite clusters with positive probability, then there is a vertex $v \in V$ such that $\P_p[L_v] > 0$.

We write $\cG_n$ for the subgraph of $\cG$ induced by $B_n$ and define the sets
\begin{align*}
  \cC_n & := \{C \subset V \mid C \text{ is a black cluster in $(\cG_n,X|_{B_n})$ such that} \\
  & \hspace{3cm} \text{there is $\{v,w\} \in E$ with $v \in C, w \in V \setminus B_n$}\}, \\
  F_n & := \bigcup_{C \in \cC_n} C ,\\
  G_n & := B_n \cap \bigcup_{C \in \cC_n} \partial^+ C ,\\
  H_n & := B_n \cap \bigcup_{C_1 \neq C_2 \in \cC_n} (\partial C_1 \cap \partial C_2).
\end{align*}
This means, that $\cC_n$ contains the black clusters of $(\cG,X)$ restricted to $B_n$ that touch the boundary of $B_n$, $F_n$ is their union, $G_n$ is the set of vertices that have a neighbour that connects to the boundary of $B_n$ via a black path and $H_n$ is the subset of $B_n$ where vertices have two neighbours in distinct black clusters connecting to the boundary of $B_n$. If a vertex $v \in B_{n-1}$ is contained in a black cluster that touches the boundary of $B_n$, then $v$ is neighbour of such a cluster if we set its color to white. This implies
\begin{align*}
  \E_p\sum_{C \in \cC_n} |C|  & = \E_p \sum_{x \in B_n} \ind\{\text{$x$ connects to the boundary of $B_n$}\} \\
  & \leq |\partial^+ B_{n-1}| + \frac{p}{1-p} \E_p \sum_{x \in B_{n-1}} \ind\{ x \in G_n\}\\
  & \leq \frac{p}{1-p} \E_p[|G_n|] + |\partial^+ B_{n-1}|
\end{align*}
and hence
\begin{align}
  \sum_{v \in B_{n-1}} \P_p[L_v] & \leq \E_p |H_n| \nonumber \\
  & \leq \E_p\bigg[  \bigg(  \sum_{C \in \cC_n} |\partial^+ C|  \bigg) - |G_n|  \bigg] \nonumber \\
  & \leq \E_p\bigg[  \sum_{C \in \cC_n} \underset{=: f(C)}{\underbrace{|\partial^+ C| - \frac{1-p}{p} |C|}}  \bigg] + \frac{1-p}{p} |\partial^+ B_{n-1}| \nonumber \\
  & = \E_p\bigg[  \sum_{C \in \cC_n} f(C)  \bigg]. \label{eq:uniqueness_via_gandolfi_grimmett_russo:1}
\end{align}
We want to apply the following large deviation result. For any $v \in V$ there is a constant $c_1(p) > 0$ that doesn't depend on $v$ such that
\begin{align}
  \P_p\left[ f(C_v) \geq \varepsilon k,\ |C_v|+|\partial^+ C_v| = k  \right] \leq e^{- c_1(p) \varepsilon^2 k}.
\end{align}
The proof can be found in \cite{gandolfi1988uniqueness2} and holds without any changes for arbitrary graphs. We define the set
\begin{align}
  \cC_n' := \{C \in \cC_n \mid |C| + |\partial^+ C| \geq \sqrt{n}\}
\end{align}
and for a fixed $\varepsilon > 0$ the event
\begin{align}
  A_n := \bigcap_{C \in \cC_n'} \{f(C) \leq \varepsilon(|C| + |\partial^+ C|)\}.
\end{align}
For the clusters that are not contained in $\cC_n'$ we have
\[
  \frac{1}{|B_n|} \E_p \sum_{C \in \cC_n \setminus \cC'_n} f(C)\leq \frac{1}{|B_n|} \E_p \sum_{C \in \cC_n \setminus \cC'_n} |\partial^+ C|\leq \frac{1}{|B_n|} \E_p \sum_{v \in B_n \setminus B_{n-\sqrt{n}}} |\partial^+ \{v\}|
\]
which tends to zero for $n \to \infty $ due to properties 1. and 2. of $\cG$.

For the clusters in $\cC_n'$ we have
\begin{align*}
  & \quad\, \frac{1}{|B_n|} \bigg( \E_p\bigg[  \sum_{C
  \in
  \cC_n'} f(C) \ \ind_{A_n} \bigg] + \E_p\bigg[   \sum_{C \in \cC_n'}
  f(C) \ \ind_{A_n^c}  \bigg] \bigg) \\
  & \leq \frac{1}{|B_n|} \bigg( \E_p\bigg[  \sum_{C
  \in \cC_n'} \varepsilon(|C| + |\partial^+ C|)  \bigg] + \E_p\bigg[
  \sum_{C \in \cC_n'} |\partial^+ C| \ \ind_{A_n^c} \bigg] \bigg) \\
  & \leq \frac{\varepsilon }{|B_n|}  \sum_{x \in
  B_n}(1+|\partial^+ \{x\}|) +
  \frac{1}{|B_n|} \sum_{x \in B_n}
  |\partial^+ \{x\}| (1 - \P_p[A_n]).
\end{align*}
The first summand tends to $\varepsilon$ for large $n$ while the second one tends to zero due to property 2. and the fact that
\begin{align*}
  1- \P_p[A_n] & \leq \P_p[\exists v \in B_n:\ C_v \in \cC_n',\
  f(C_v) > \varepsilon(|C_v| + |\partial^+ C_v|)] \\
  & \leq \sum_{v \in B_n} \sum_{k \geq \sqrt{n}} \P_p[f(C_v) >
  \varepsilon k,\ |C_v| + |\partial^+ C_v| = k] \\
  & \leq |B_n| \sum_{k \geq \sqrt{n}} e^{-c_6 \varepsilon^2 k}.
\end{align*}
Putting everything together we see, that
\[
  \lim_{n\to \infty } \frac{1}{|B_n|} \sum_{v \in B_n} \P_p[L_v] = 0
\]
and assumption 1.\ implies the assertion.\qed
\bigskip

The assumptions in Theorem \ref{thm:uniqueness_GGR} resemble some kind of ergodicity on $\cG$ and may in fact be shown for graphs induced by ergodic random tessellations where balls have an a.s.\ polynomial growth, i.e.\ if there are constants $c_4, c_5 > 0$ such that $\lim_{n\to \infty } |B_n| n^{-c_4} = c_5$. However it is not trivial to show this polynomial growth for an arbitrary random tessellation. We will address this problem in our second paper on first passage percolation on random tessellations. Note also that property 2. doesn't depend on the choice of the root $\mathbf{0}$.

\subsection*{The planar case}
The third result we want to generalize in this section is, that in a planar lattice there can't be a coexistence of an infinite white and an infinite black cluster. The most basic proof of this result can be found in \cite[p. 289]{grimmett1999percolation} named argument of Zhang. It was later generalized to lattices with a $k$-fold symmetry in \cite{bollobas2008kfoldsymmetry}. We will use some of their arguments, to show, that the uniqueness of the infinite cluster already implies that $p_c \geq \tfrac{1}{2}$ in the planar case.
\begin{Thm}\label{thm:uniqueness_implies_p_c_geq_half}
  Let $\cG = (V, E)$ be an infinite connected planar graph. There is at most one infinite cluster at $p = \tfrac{1}{2}$ if and only if $\theta_{\frac{1}{2}}(v) = 0$ for any $v \in V$. Moreover we have
  \[
    p_c \geq \frac{1}{2}
  \]
  in this case.
\end{Thm}

The proof of Theorem \ref{thm:uniqueness_implies_p_c_geq_half} contains some tedious topological details. A background in planar graph theory with a quite rigorous approach can be found in \cite{diestel2010graph} or \cite{mohar2001graphs}. The Jordan Curve Theorem (JCT) can be found in \cite{hales2007jordan} along with some interesting historical remarks.

We call a continuous map from $[0,1] \to \R^d$ a \emph{curve}. A curve $\varphi$ is called \emph{closed} if $\varphi(0) = \varphi(1)$, \emph{Jordan curve} if it is injective, \emph{closed Jordan curve} if it is closed and injective when restricted to $[0,1)$ and \emph{polygonal} if it is piecewise linear. We will identify any curve with its image in $\R^d$. The JCT states that for any closed Jordan curve $\varphi$ the set $\R^2 \setminus \varphi$ consists of one bounded (the interior) and one unbounded (the exterior) connected component. Moreover if two points $x, y \in \R^2$ are connected by a Jordan curve that crosses $\varphi$ an odd number of times, then one of these points lies in the interior and one lies in the exterior of $\varphi$.

A \emph{cut vertex} of a connected graph $\cG$ is a vertex $v$ such that deleting $v$ results in $\cG$ being not connected anymore. A graph $\cG$ is called planar if there is an embedding of $\cG$ in the plane such that all edges are piecewise linear and don't intersect (edges do not contain their endpoints). Moreover, the embedding has to be locally finite, i.e.\ any bounded component of $\R^2$ is intersected only by a finite number of vertices and line-segments of the embedding.

We will state two Lemmas first, that contain the topological arguments.

\begin{Lemma}\label{le:jordan-curve_around_B_n}
  Let $\cG = (V, E)$ be an infinite connected planar graph with root $\mathbf{0}$. Then for all $n \in \N$ there is a closed Jordan curve $\varphi := \varphi(n, \cG)$ with the following properties:
  \begin{enumerate}
    \item No edge of $\cG$ intersects $\varphi$.
    \item Each vertex in $B_n$ is either contained in $\varphi$ or in its interior.
    \item We have $R_n \subset \varphi$, where $R_n$ is the set of vertices of $B_n$ where an infinite path may start, that intersects $B_n$ only once.
  \end{enumerate}
\end{Lemma}
\emph{Proof:} We will consider a finite connected planar graph $\cG_f = (V_f, E_f)$ first. A \emph{closed walk} is a cycle in $\cG_f$ that allows to visit vertices multiple times. Any face of $\cG_f$ induces a walk along its boundary in a natural way. This is called the \emph{facial walk}. It is well known (see \cite{diestel2010graph}) that $\cG_f$ has exactly one unbounded face $F$ with a facial walk $L = (l_1, \dots, l_m)$, $m \in \N$. The walk $L$ has the property that two consecutive edges $e_i = \{l_i, l_{i+1}\}$ and $e_{i + 1} = \{l_{i+1}, l_{i+2}\}$, $i \in [m-2]$ enter $l_{i+1}$ in clockwise order without any other edges in between. Moreover, there is a small circular sector enclosed by the ends of the edges $e_i$ and $e_{i+1}$ that is contained in $F$. If we fix a starting vertex and require each edge to be traversed at most once, the walk $L$ is uniquely determined.
\begin{figure}
  \centering
  \includegraphics{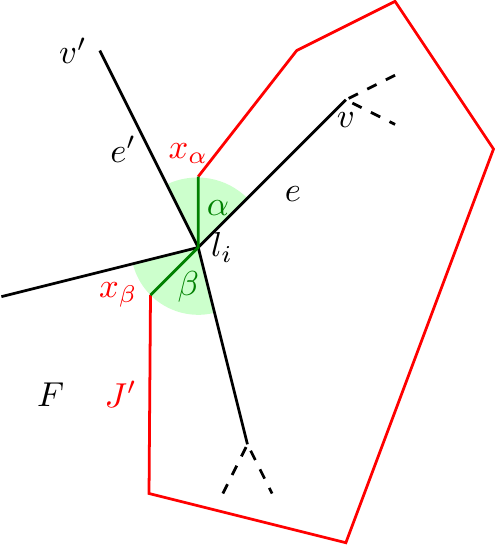}
  \includegraphics{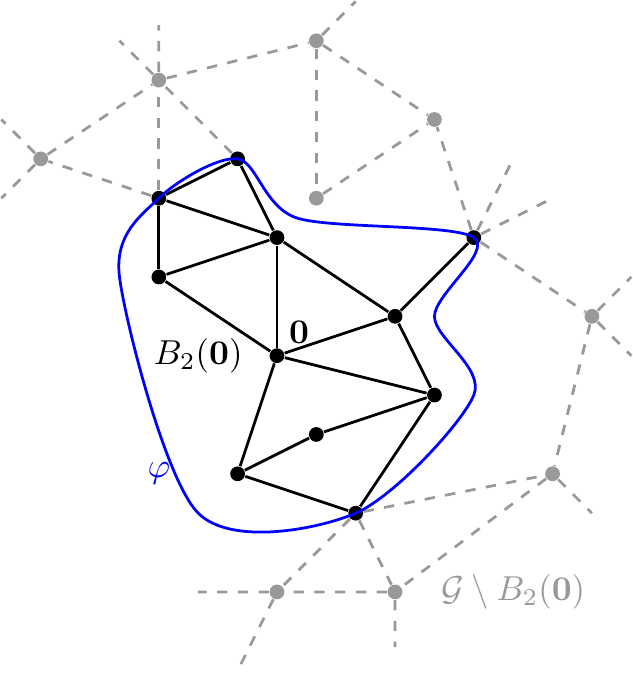}
  \caption{The left figure shows the situation in the first claim. The right figure shows an example of the whole situation.}
  \label{fig:cut-vertex_in_facial_walk_and_jordan_curve_around_B_n}
\end{figure}

Our first claim is, that any vertex that is contained multiple times in $L$ is a cut-vertex of $\cG_f$. Assume $l_i$ is contained at least two times in $L$, then there are two circular sectors $\alpha, \beta$ at $l_i$ that are contained in $F$. We connect $\alpha$ and $\beta$ with a polygonal curve $J'$ starting in $x_{\alpha}$, ending in $x_\beta$ that is contained in $F$ (see Figure \ref{fig:cut-vertex_in_facial_walk_and_jordan_curve_around_B_n}). Extending $J'$ by line-segments from $l_i$ to $x_\alpha$ and from $x_\beta$ to $l_i$ we obtain a closed Jordan curve $J$.

Considering the two edges $e := \{l_i, v\}$ and $e' := \{v', l_i\}$ that enclose $\alpha$ we observe that the curve starting in $v$ traversing $e$ as well as $e'$ and ending in $v'$ crosses $J$ exactly once at $l_i$. Hence $v$ and $v'$ lie in different connected components of $\R^2 \setminus J$ and the deletion of $l_i$ implies that $v$ and $v'$ are not connected anymore.

Now we look at the finite subgraph $\cG_n$ of $\cG$ that is induced by $B_n$. The the facial walk $L_n$ of the unbounded component $F_n$ of $\cG_n$ is the natural basis for the closed curve $\varphi$. We will construct $\varphi$ in three steps. We start with the closed curve $\varphi_1$ which is obtained by traversing $L_n$ along its edges once. We obtain $\varphi_2$ by replacing each part of $\varphi_1$ going from a vertex $l_i$ to $l_{i+1}$ with a curve that also connects $l_i$ and $l_{i+1}$ but lies in $F_n$ and doesn't intersect $\cG$ apart from the starting and the ending points. Due to the local finiteness of the embedding of $\cG$ this is even possible in such a way, that the segments in $\varphi_2$ from $l_i$ to $l_{i+1}$ doesn't intersect for different $i$. The curve $\varphi_2$ already fulfills properties 1., 2. and 3. as $R_n$ is a subset of the boundary of $F_n$. However it is still possible, that a vertex occurs multiple times in $L_n$. In this case $\varphi_2$ is not injective.
\begin{figure}
  \centering
  \includegraphics{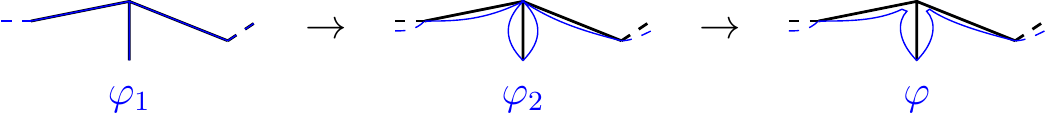}
  \caption{How $\varphi_1$, $\varphi_2$ and the final $\varphi$ might look like.}
  \label{fig:evolution_of_phi}
\end{figure}

To correct this problem, we modify $\varphi_2$ to get $\varphi$ by skipping any points $l_i$ that occur more than once. To skip a point $l_i$, we delete a small part of $\varphi_2$ that leads to and a small part of $\varphi_2$ that comes from $l_i$ and connect the dangling ends directly without visiting $l_i$ and without intersecting anything else (see Figure \ref{fig:evolution_of_phi}). It remains to show, that $R_n$ contains no cut-vertex, as in this case we don't destroy property 3. by the procedure.

Let $v$ be a cut-vertex of $G_n$ and let $S_1,\dots, S_k$, $k \geq 2$ be the connected components of $G_n$ that emerge after the deletion of $v$. Let $i \in [k]$ be such that $\mathbf{0} \notin S_i$. Hence for any $w \in S_i$ we have $n \geq d_{\cG_n}(w,\mathbf{0}) > d_{\cG_n}(v,\mathbf{0})$. However for any vertex $u \in R_n$ we know that $n = d_{\cG_n}(u,\mathbf{0})$. This implies that $R_n$ contains no cut-vertex of $\cG_n$ which finishes the proof.\qed
\bigskip

A \emph{Jordan ray} $\gamma$ is an injective mapping $\gamma:[0,\infty ) \to \R^2$ such that $\lim_{t\to \infty } \|\gamma(t)\|_2 = \infty $. We identify these objects again with their image in $\R^2$.

\begin{Lemma}\label{le:diagonals_cross}
  Let $\varphi$ be a closed Jordan curve containing the points $v_1, \dots, v_4$ in this order. Let $\gamma_1, \dots, \gamma_4$ be polygonal Jordan rays and let $\rho_{13}, \rho_{24}$ be two polygonal Jordan curves with the following properties:
  \begin{enumerate}
    \item $(\gamma_1 \cup \gamma_3) \cap (\gamma_2 \cup \gamma_4) = \emptyset $,
    \item $\gamma_i \cap \varphi = \{\gamma_i(0)\} = \{v_i\}$ für $i \in [4]$,
    \item $\rho_{13}(0) \in \gamma_1$, $\rho_{13}(1) \in \gamma_3$,
    \item $\rho_{24}(0) \in \gamma_2$, $\rho_{24}(1) \in \gamma_4$.
  \end{enumerate}
  Then we have
  \[
    (\gamma_1 \cup \rho_{13} \cup \gamma_3) \cap (\gamma_2 \cup \rho_{24} \cup \gamma_4) \neq \emptyset .
  \]
\end{Lemma}
\emph{Proof:} We chose a radius $r$ large enough such that $\rho_{13}, \rho_{24}, \varphi \subset [-r+1,r-1]^2$. Let $w_i := \gamma_i(t_i)$ where $t_i := \min\{t \in [0, \infty ) \colon \|\gamma_i(t)\|_\infty = r\}$, $i \in [4]$ be the first intersection point of $\gamma_i$ and $\partial [-r,r]^2$. It follows with the help of the JCT that we may traverse $\partial [-r,r]^2$ in a way, such that $w_1, \dots, w_4$ are visited in this order. The sets $\gamma_1 \cup \rho_{13} \cup \gamma_3$ and $\gamma_2 \cup \rho_{24} \cup \gamma_4$ contain by construction polygonal curves $\rho_{13}'$ starting in $w_1$, ending in $w_3$, being contained in $[-r, r]^2$ and $\rho_{24}'$ starting in $w_2$, ending in $w_4$, being contained in $[-r, r]^2$. Another application of the JCT yields, that $\rho_{13}'$ and $\rho_{24}'$ have to intersect, due to the order in which $w_1, \dots, w_4$ lie on $\partial [-r,r]^2$.\qed
\bigskip
\begin{figure}
  \centering
  \includegraphics{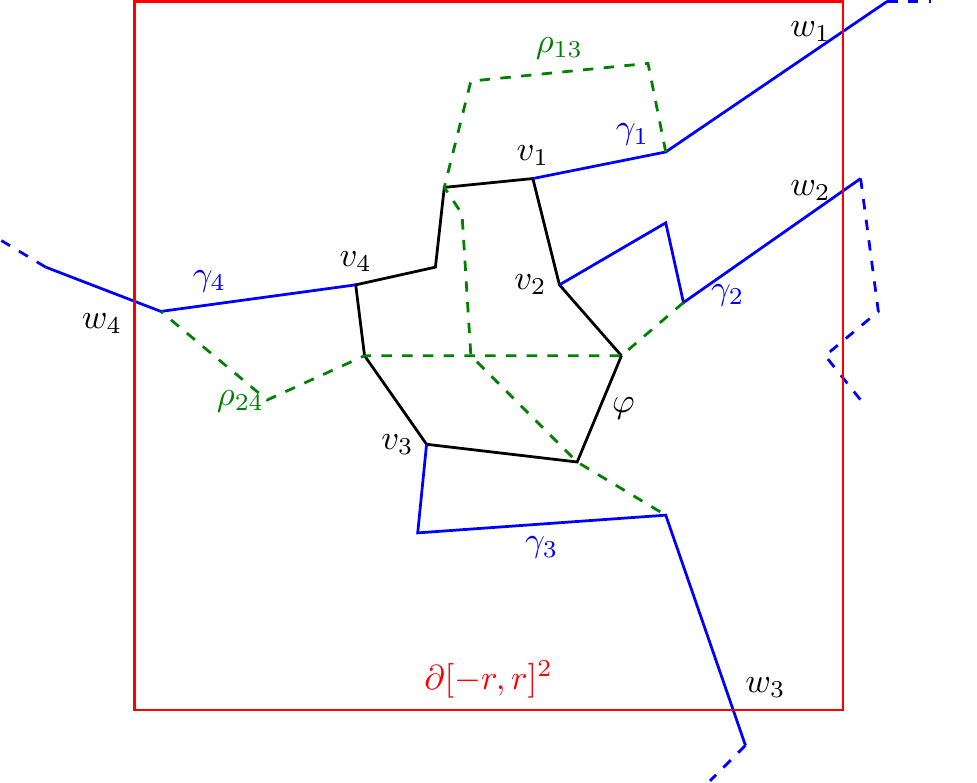}
  \caption{Setup of Lemma \ref{le:diagonals_cross}.}
  \label{fig:diagonals_cross}
\end{figure}

\emph{Proof of Theorem \ref{thm:uniqueness_implies_p_c_geq_half}:} If $\theta_{\frac{1}{2}}(v) = 0$ for all $v$, then there is a.s.\ no infinite cluster which proves one direction.

Let $p = \tfrac{1}{2}$, $\varepsilon > 0$ and let us assume that there is a.s.\ exactly one infinite black cluster $C_{\infty }^+$. This implies the a.s.\ existence of an infinite white cluster $C_\infty ^-$. We may choose $n \in \N$ large enough such that
\begin{align}
  \P_{\frac{1}{2}}[R_n \cap C_\infty^+ \neq \emptyset ] = \P_{\frac{1}{2}}[B_n \cap C_\infty^+ \neq \emptyset ] > 1 - \varepsilon,
\end{align}
where $R_n$ is defined as in Lemma \ref{le:jordan-curve_around_B_n}. Let $\varphi := \varphi(n, \cG)$ be the closed Jordan curve that exists due to Lemma \ref{le:jordan-curve_around_B_n} for the chosen $n$. Let $(r_1, \dots, r_m)$, $m \in \N$ be the points of $R_n$ ordered in the way that is induced if $\varphi$ is traversed in clockwise direction.

Now a slight adaptation of the arguments from \cite{bollobas2008kfoldsymmetry} is enough to finish the proof. We state it for the convenience of the reader. For any set $U \subset R_n$ we define the decreasing event $E^+(U) :=$"There is no infinite black path starting in $U$ that intersects $R_n$ only once." and the function
\begin{align}
  f(U) := \P_{\frac{1}{2}}[E^+(U)].
\end{align}
Clearly the following properties hold for $f$:
\begin{enumerate}[(i)]
  \item $f(\emptyset ) = 1$,
  \item $f(\{v\}) \geq \P_{\tfrac{1}{2}}[X_v = 0] = \frac{1}{2}$, $v \in R_n$,
  \item $f(U \cup W) \geq f(U) f(W)$, $U, W \subset R_n$ (due to FKG-inequality),
  \item $f(U) \geq f(W)$, $U \subset W\subset R_n$,
  \item $f(R_n) \leq \varepsilon$.
\end{enumerate}
It follows from (ii), (iii) and (iv) that for any $U \subset R_n$ and $v \in R_n$
\begin{align}\label{eq:uniqueness_implies_p_c_geq_half:1}
  f(U) \geq f(U \cup \{v\}) \geq \frac{1}{2} f(U).
\end{align}
The idea is, to separate $R_n$ into four parts $T_1, \dots, T_4$ that have a low $f$-value. At first we separate $R_n$ into two parts $T_{12}^{(k)} := (r_1,\dots, r_k)$ and $T_{34}^{(k)} := (r_{k+1}, \dots, r_m)$, $k \in [m]$. Taking $k = m$ implies $f(T_{12}^{(k)}) = f(R_n) \leq \varepsilon$ and $f(T_{34}^{(k)}) = f(\emptyset ) = 1$. It follows from \eqref{eq:uniqueness_implies_p_c_geq_half:1}, that if we reduce $k$ by one, we increase the value of $f(T_{12}^{(k)})$ by a factor of at most two while we decrease $f(T_{34}^{(k)})$ by a factor of at most two. Hence there has to be a $k \in [m]$ such that
\begin{align}
  f(T_{12}^{(k)}) \leq f(T_{34}^{(k)}) \leq 4 f(T_{12}^{(k)}).
\end{align}
For such a $k$
\begin{align*}
  f(T_{12}^{(k)})^2 &\leq f(T_{12}^{(k)}) f(T_{34}^{(k)}) \leq f(R_n) \leq \varepsilon, \\
  f(T_{34}^{(k)})^2 &\leq 4 f(T_{12}^{(k)}) f(T_{34}^{(k)}) \leq 4 f(R_n) \leq 4 \varepsilon.
\end{align*}
In the same way, we separate $T_{12}^{(k)}$ and $T_{34}^{(k)}$ another time into the parts $T_1, T_2$ and $T_3, T_4$ such that
\[
  f(T_i) = \P_{\frac{1}{2}}[E^+(T_i)] = \P_{\frac{1}{2}}[E^-(T_i)] \leq 2 \sqrt{2} \sqrt[4]{\varepsilon}, \quad i \in [4],
\]
where we define $E^-(U)$ for $U \subset R_n$ in same way as $E^+(U)$ except that we use white paths instead of black ones. Due to symmetry we have $\P_{\frac{1}{2}}[E^-(U)] = \P_{\frac{1}{2}}[E^+(U)]$ and by choosing $\varepsilon$ small enough we obtain, that the event
\[
  E^+(T_1)^c \cap E^-(T_2)^c \cap E^+(T_3)^c \cap E^-(T_4)^c
\]
has positive probability. This event describes the case where there are black infinite paths $\gamma_1, \gamma_3$ emanating from $T_1$ and $T_3$ and white infinite paths $\gamma_2, \gamma_4$ emanating from $T_2$ and $T_4$. Each of these paths intersects $B_n$ only in its starting point. We assumed that there is only one infinite cluster of each color and hence $\gamma_1$ and $\gamma_3$ has to be connected by some finite black path $\rho_{13}$, while same is true for some white finite path that has to connect $\gamma_2$ and $\gamma_4$. These paths fulfill by construction exactly the required properties in Lemma \ref{le:diagonals_cross} and hence $\gamma_1 \cup \rho_{13} \cup \gamma_3$ and $\gamma_2 \cup \rho_{24} \cup \gamma_4$ have a nonempty intersection, which is a contradiction as it would imply a vertex that has both colors.\qed
\bigskip

\section{Percolation on random tessellations}\label{sec:perc_on_random_tess}
In this section we want to define our model in a rigorous way, discuss several reasonable assumptions and show the uniqueness of the infinite cluster by adapting the Burton Keane argument from \cite{burton1989density}.

We will interpret randomly colored random tessellations as independently marked particle process. This leads to the following definitions and notational conventions (a broad introduction to point processes and random tessellations can be found in \cite{schneider2008stochastic}).

Let $D$ be a metric space equipped with the Borel-$\sigma$-algebra $\cB(D)$. We write $\mathbf{N}(D)$ for the set of locally finite counting measures on $D$ and equip it with the $\sigma$-algebra $\cN(D)$ generated by the sets $\{\eta \in \mathbf{N}(D) \mid \eta(A) = k\}$, $A \in \cB(D), k \in \N_0 := \N \cup \{0\}$. A measure $\eta \in \mathbf{N}(D)$ is called locally finite iff $\eta(A) < \infty $ for all bounded $A \in \cB(D)$. A measurable mapping $\Phi:\Omega \to \mathbf{N}(D)$ is called a \emph{point process} on $D$ and is to be interpreted as a random collection of points in $D$. Each point process permits a representation
\[
  \Phi = \sum_{i = 1}^{\Phi(D)} \delta_{\zeta_i},
\]
where $\delta$ is the Dirac measure and $(\zeta_i)_{i \in \N}$ are $D$-valued random variables \cite[Lemma 3.1.3]{schneider2008stochastic}. The measure $\Theta := \E\Phi$ on $D$ is called the \emph{intensity measure} of $\Phi$. In the important special case, that $D = \R^d$ and $\Theta = \gamma \lambda^d$ we call $\gamma$ the \emph{intensity} of $\Phi$ ($\lambda^d$ is the Lebesgue measure). Random tessellations will later be defined by letting $D$ be the space $\cC^d$ of compact and convex subsets of $\R^d$ equipped with the Hausdorff metric.

To be able to add a color information to each cell, it is convenient to work with marked point processes. If we have a point process $\Phi = \{\zeta_1, \zeta_2, \dots\}$ on $D$ with representation $\sum_{i = 1}^{\Phi(D)} \delta_{\zeta_i}$ and an i.i.d. sequence $X = (X_i)_{i \in \N}$ with $X_1 \sim \Ber(p)$, $p \in [0,1]$ that is independent of $\Phi$, we call
\[
  \Phi_X := \sum_{i = 1}^{\Phi(D)} \delta_{(\zeta_i, X_i)}
\]
the \emph{independently marked version} of $\Phi$. To indicate which $p$ is used, we will write $\P_{p}$ instead of $\P$ where necessary.

From now on we will only work with point processes on $\R^d$ or $\cC^d$ and their independently marked versions (each marked point process on $D$ is also a point process on $D \times \{0,1\}$). If $D$ is equal to either of these spaces, the canonical translation operator $T_x:\mathbf{N}(D) \to \mathbf{N}(D)$, $x \in \R^d$ is defined by
\[
  T_x \eta(A) := \eta(A - x), \quad \eta \in \mathbf{N}(D), A \in \cB(D).
\]
We use the same notation for the shift $T_x:\mathbf{N}(D \times \{0,1\}) \to \mathbf{N}(D \times \{0,1\})$ on the marked spaces defined by
\[
  T_x \eta(A \times B) := \eta( (A - x) \times B), \quad \eta \in \mathbf{N}(D), A \in \cB(D), B \subset \{0,1\}.
\]
This corresponds to the idea that only the points are shifted while each point retains its mark. A point process on $\cC^d$ is also called a \emph{particle process}.

We recall, that a point process $\Phi$ is \emph{stationary} iff $T_x \Phi \overset{d}{=} \Phi$ for all $x \in \R^d$. Let $\cI$ be the $\sigma$-Algebra of translation invariant events, i.e.\ events $A \in \cN(D)$ with $T_x A = A$ for all $x \in \R^d$. A stationary point process $\Phi$ is called \emph{ergodic} if $\P[\Phi \in A] \in \{0,1\}$ for all $A \in \cI$. It can be shown, that if $\Phi$ is ergodic, then its independently marked version $\Phi_X$ is ergodic too; see \cite[Proposition 12.3.VI.]{daley2007introductionvolzwei}.

We recall the definition of the \emph{Laplace functional} of a point process $\Phi$ applied to a function $f:D \to [0,\infty )$
\[
  L_{\Phi}(f) := \E \exp\Big( -\int f(x)\ \Phi(dx)\Big).
\]
We will extend this definition to $\R$-valued functions $f$ and remark that the integral or the expectation might not exist in this case. An introduction to point processes and the Laplace functional for random measures can be found in \cite{daley2003introduction}.

After the introduction of point and particle processes, we turn to tessellations. A set $Z \in \cC^d$ with non-empty interior is called a cell. A countable set $m := \{Z_1, Z_2, \dots\}$ of cells is called a \emph{tessellation} (or \textbf{m}osaic) if
\begin{enumerate}
  \item each ball in $\R^d$ is intersected by at most a finite number of cells of $m$,
  \item the cells of $m$ cover $\R^d$,
  \item the interiors of any two distinct cells in $m$ doesn't overlap.
\end{enumerate}
The cell of $m$ that contains $x \in \R^d$ is denoted by $Z_x(m)$ (if there is more than one cell containing $x$ we chose an arbitrary rule to break ties). The cell $Z_0(m)$ is called the \emph{zero cell}.
\begin{figure}
  \centering
    \includegraphics[scale=.9]{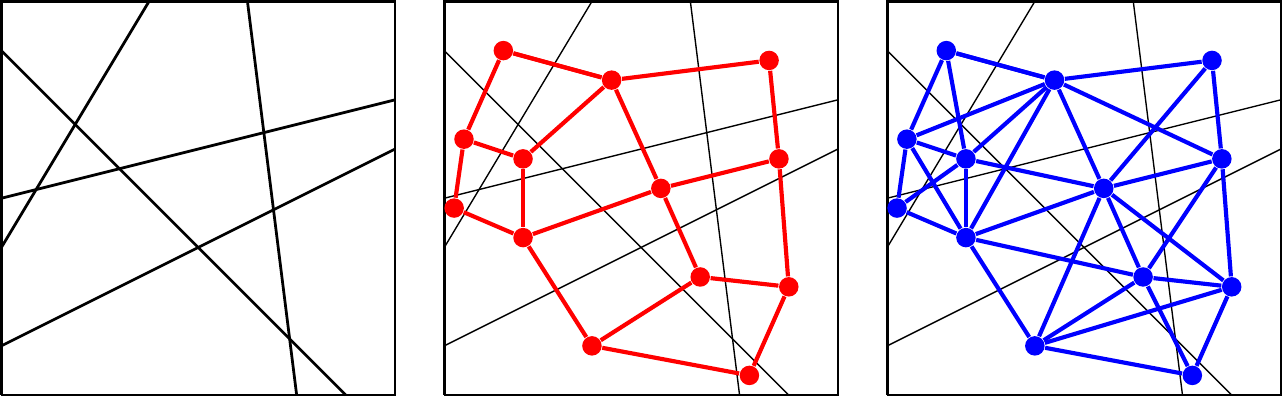}
    \caption{A section of a Poisson hyperplane tessellation $m$ and its induced graphs \textcolor{red}{$\cG_{\mathfrak{m}}$} und \textcolor{blue}{$\cG_{\mathfrak{m}}^*$}}
  \label{fig:tessellation_and_induced_graphs}
\end{figure}

Let $\mathbf{M} \subset \mathbf{N}(\cC^d)$ be the set of tessellations and observe, that any tessellation $m \in \mathbf{M}$ induces a graph $\cG_m := (m, E_m)$ with vertex set $m$. Two cells $Z_1, Z_2 \in m$ are adjacent in $\cG_m$ iff they have a $(d-1)$-dimensional intersection, i.e.\ iff $Z_1 \cap Z_2$ is not contained in any $(d-2)$-dimensional hyperplane. The tessellation $m$ induces a second graph $\cG_m^* := (m, E_m^*)$ where any two cells with nonempty intersection are adjacent. The distinction of $\cG_m$ and $\cG_m^*$ will mostly be relevant in the 2-dimensional case. Apart from that, all results will hold for both graphs. The zero cell $Z_0(m)$ is the root $\mathbf{0}$ in $\cG_m$ and $\cG_m^*$.

We denote by $\mathfrak{F}(\mathbf{M}) := \cN(\cC^d)|_\mathbf{M}$ the trace of $\cN(\cC^d)$ on $\mathbf{M}$ and call a measurable mapping $M:\Omega \to \mathbf{M}$ a random tessellation. Hence a random tessellation is a point process of convex compact particles that form a tessellation. In the same spirit let $\mathbf{M}_c := \{(Z_i, X_i)_{i \in \N} \mid (Z_i)_{i \in \N} \in \mathbf{M},\ X_j \in \{0,1\},\ j \in \N\}$ be the set of \emph{colored tessellations} with the $\sigma$-Algebra $\mathfrak{F}(\mathbf{M}_c) := \cN(\cC^d \times \{0,1\})|_{\mathbf{M}_c}$. For $m = \{Z_1, Z_2, \dots\} \in \mathbf{M}$ and a random or deterministic $\{0,1\}$-valued sequence $X = \{X_1, X_2, \dots\}$ we define the marked tessellation $m_X := \{(Z_1, X_1), (Z_2, X_2), \dots\} \in \mathbf{M}_c$. A marked tessellation $m_X$ induces a colored graph $\cG_{m,X} := (\cG_{m}, c)$ where $c(Z_i) := X_i$, $i \in \N$ (the same notations are used for a random tessellations $M$ in place of $m$).
\begin{figure}
  \centering
    \includegraphics[scale=.9]{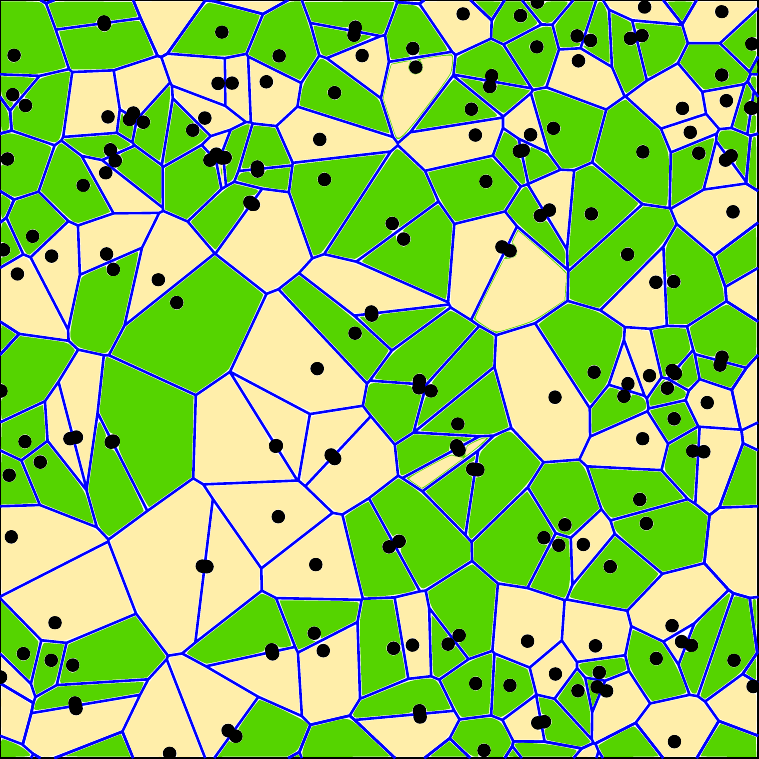}
    \caption{A section of a randomly colored ($p = 0.6$) Voronoi-tessellation that was created from a Poisson cluster process (Data by Michael Klatt).}
  \label{fig:percolation_example_on_voronoi-tessellation}
\end{figure}

For the rest of this article we will work with a stationary random tessellation $M = \{Z_1, Z_2, \dots\}$ and a random sequence $X = (X_1, X_2, \dots)$. Under $\P_p$ the random tessellation $M$ is independent of the i.i.d. sequence $X$ that has $\Ber(p)$, $p \in [0,1]$ distributed marginals. Hence $M_X$ is stationary and we have
\begin{align}\label{eq:doubly-stochastic-model}
  \P_{p}[M_X \in A] = \int_{\mathbf{M}} \P_{p}[m_X \in A]\ \P_{M}(dm), \quad A \in \mathfrak{F}(\mathbf{M}_c)
\end{align}
where $\P_{M}$ denotes the distribution of $M$.

The first thing one might observe is, that unlike in the case of percolation on a fixed graph, in our model the existence of an infinite black cluster can have a probability different from 0 or 1. We could, for instance, create a random tessellation by taking a randomly shifted square lattice with probability 1/2 and a randomly shifted honeycomb lattice otherwise. The resulting random tessellation would exhibit an infinite black cluster with probability 1/2 for $p$ between the lattice dependent percolation thresholds (which are known to be different). To rule out this somehow pathological case, we will restrict ourselves to ergodic random tessellations most of the time. As the set $A \subset \mathbf{M}_c$ of colored tessellations that contain an infinite black cluster is translation invariant, we have $\P_{p}[M_X \in A] \in \{0,1\}$ for any ergodic random tessellation $M$, since $M_X$ is also ergodic in this case.

Now it is only natural to define the percolation function and threshold for $M$ by
\[
  \theta_p(M) := \P_{p}[|C(\mathbf{0}, \cG_{M,X})| = \infty ]
\]
and
\[
  p_c(M) := \sup\{p \in [0,1] \mid \theta_p(M) = 0\}.
\]
By a standard coupling argument we see that $\theta_p(M)$ is non-decreasing in $p$. It is also clear that there is an infinite black cluster in $\cG_{M,X}$ for any $p > p_c(M)$ a.s.\ and that there is no infinite black cluster in $\cG_{M,X}$ for $p < p_c(M)$ a.s.\ . The existence of an infinite cluster at $p_c(M)$ is obviously an open and hard problem in most cases.

\subsection*{Uniqueness of the infinite cluster with Burton and Keane}
We are now in the position to adapt the Burton Keane argument \cite{burton1989density} to our model.

\begin{Thm}\label{thm:uniqueness-burton-keane}
  Let $M$ be a stationary random tessellation and $p \in [0,1]$. If
  \begin{align}\label{eq:burton-keane:moment-assumption}
    \E_{p}[|\{Z \in M \mid Z \cap [0,1]^d \neq \emptyset \}|] < \infty
  \end{align}
   then there is a.s.\ at most one infinite black cluster in $\cG_{M,X}$.
\end{Thm}
\emph{Proof:} The claim is trivial for $p \in \{0,1\}$, so let $p \in (0,1)$ for the rest of the proof. First, we assume that $M$ is ergodic. We already mentioned that this implies the ergodicity of $M_X$. For each $n \in \N_0 \cup \{\infty \}$ we define the set
\[
  E_n := \{m_c \in \mathbf{M}_c \mid \text{there are exactly $n$ infinite black clusters in $\cG_{m,c}$}\}.
\]
These sets are translation invariant and hence exactly one of the events $\{M_X \in E_n\}$, $n \in \N_0 \cup \{\infty \}$ will a.s.\ hold while all others will a.s.\ not.

Let us assume that $\{M_X \in E_n\}$ holds a.s.\ for a fixed $2 \leq n < \infty $. In this case $M$ lies a.s.\ in the set
\[
  A := \{m \in \mathbf{M} \mid \P_{p}[\text{there are $n$ infinite black cluster in $\cG_{m,X}$}] = 1\}.
\]
We fix an $m \in A$ and define the random variables $N(i, r)$, $i \in \{0,1\}$, $r \in \N$ as the number of infinite black clusters in the colored graph which we obtain, if all vertices from $\cG_{m,X}$ in $B_r(\mathbf{0}, \cG_m)$ are set to the value $i$. The probability, that all vertices in $B_r(\mathbf{0}, \cG_m)$ have the value $i$ in $\cG_{m,X}$ is positive and hence
\[
  \P_{p}[N(0,r) = N(1,r) = n] = 1, \quad r \in \N.
\]
As $n < \infty $ the random variables $N(0,r)$ and $N(1,r)$ can only be equal, if $B_r(\mathbf{0}, \cG_m)$ is intersected by at most one infinite black cluster of $\cG_{m,X}$. However the number of clusters intersecting $B_r$ converges a.s.\ to $n$ if $r$ tends to infinity, which leads to a contradiction.

It remains to rule out the case of an infinite number of infinite black clusters. Let $m_c \in \mathbf{M}_c$ and $x \in \R^d$. We call $x$ a \emph{trifurcation point of $m_c$ with parameters $(r_1, r_2)$} if
\begin{itemize}
  \item all vertices in $B_{r_1}(Z_x(m),\cG_m)$ are black in $\cG_{m,c}$,
  \item the outer boundary $\partial^+ B_{r_1}(Z_x(m), \cG_m)$ is intersected by at least three infinite black clusters in the colored graph we obtain from $\cG_{m,c}$ by setting all values of vertices in $B_{r_1}(Z_x(m), \cG_m)$ to 0,
  \item all cells in $B_{r_1}(Z_x(m), \cG_m)$ are contained in $[-r_2, r_2]^d$.
\end{itemize}
If $M \in E_\infty $ a.s.\ then we might choose $r_1$ large enough such that with positive probability $B_{r_1}(Z_0(m), \cG_m)$ is intersected by at least three infinite black clusters. Hence by choosing $r_1$ and $r_2$ large enough we can ensure that each $x \in \R^d$ has the same positive probability to be a trifurcation point of $M_X$ with parameters $(r_1, r_2)$.

If $x \in \R^d$ is a trifurcation point of $m_c \in \mathbf{M}_c$ with parameters $(r_1, r_2)$, then setting the color of all vertices in $B_{r_1}(Z_x(m),\cG_m)$ to white splits the former black cluster of $Z_x$ into at least three infinite black clusters. By applying a standard combinatorial lemma (see e.g. \cite[p. 121, Lemma 3]{bollobas2006percolation}) we conclude that if there are at least $k \in \N$ points of the set $3r_2 \Z^d \cap [-r,r]^d$, $r > 0$ trifurcation points of $m_c$ with parameters $(r_1, r_2)$, then there are at least $k+2$ disjoint infinite black paths in $\cG_{m,c}$ starting in a cell contained in $[-r -r_2, r + r_2]$. In this case the boundary of $[-r -r_2, r + r_2]$ has to be intersected by at least $k$ cells.

This leads again to a contradiction if $r$ tends to infinity as the expected number of trifurcation points of $M_X$ with parameters $(r_1, r_2)$ that are contained in $3r_2 \Z^d \cap [-r,r]^d$ is of order $r^d$ while the expected number of cells that intersect the boundary of $[-r -r_2, r + r_2]$ is at most of order $r^{d-1}$.

To lift this result to the case where $M$ is stationary, we recall that there is an appropriate $\sigma$-algebra and a measure $\nu$ on the set $\cL$ of all distributions of ergodic random tessellations such that
\begin{align}\label{eq:burton-keane:ergodic-decomposition}
  \P[M \in A] = \int_{\cL} P(A)\ \nu(dP), \quad A \in \mathfrak{\mathbf{M}}.
\end{align}
It follows from our moment assumption \eqref{eq:burton-keane:moment-assumption}, that for $\nu$-a.e. $P \in \cL$
\[
  \int_{\mathbf{M}} |\{Z \in m \mid Z \cap [0,1]^d \neq \emptyset \}|\ P(dm) < \infty .
\]
Hence applying \eqref{eq:burton-keane:ergodic-decomposition} to the set of tessellations $m$ where $\cG_{m,X}$ $\P_{p}$-a.s.\ contains at most one infinite black cluster proves the assertion.\qed\bigskip

This result is immediately applicable to the two dimensional case, where we obtain the following corollary.

\begin{Cor}
  Let $M$ be a stationary random tessellation of $\R^2$. If $\E[|\{Z \in M \mid Z \cap [0,1]^2 \neq \emptyset \}|] < \infty $, then $\theta_{\frac{1}{2}}(M) = 0$ and
  \[
    p_c(M) \geq \frac{1}{2}
  \]
\end{Cor}
\emph{Proof:} The graph $\cG_M$ is a.s.\ planar and connected. Due to Theorem \ref{thm:uniqueness-burton-keane} for a.e. $m \in \mathbf{M}$ the colored graph $\cG_{m,X}$ contains $\P_{p}$-a.s.\ at most one infinite black cluster. Hence by applying Theorem \ref{thm:uniqueness_implies_p_c_geq_half} to $\cG_m$ we have that there is $\P_{\frac{1}{2}}$-a.s.\ no infinite cluster in $\cG_{m,X}$ for a.e. $m \in \mathbf{M}$.\qed\bigskip

One might guess, that for many ergodic random tessellations $M$ of $\R^2$ with $\cG_M = \cG_M^*$ a.s.\ the critical value is exactly one half, but this has been shown only for the Poisson-Voronoi tessellation in two dimensions \cite{bollobas2006critical}. We consider it an interesting open problem to generalize this result to other 2-dimensional tessellations.

At the end of this Section, we want to give a counterexample, that shows that ergodicity and even mixing is by itself not enough to ensure a non-trivial phase transition. We recall, that a stationary random tessellation $M$ is called \emph{mixing}, if
\begin{align}\label{eq:mixing}
  \lim_{\|x\|\to \infty } \P[M \in A \cap T_x B] = \P[M \in A] \P[M \in B]
\end{align}
for any $A, B \in \mathfrak{F}(\mathbf{M})$.

In \cite{heil2012stationary} Heil constructs a random partition of $\Z^2$ called streetgrid, that is mixing (with respect to $(T_x)_{x \in \Z^2}$). Each element of the partition is a rectangle. Hence we could think of this partition as a random mixing tessellation (though one would have to circumvent the discreteness is his construction, to be rigorous). The construction has the property that in a.e. realization there are an infinite number of disjoint quadruples of rectangles that enclose the origin. By enclose we mean, that any path starting from the origin intersects at least one rectangle of each quadruple. If we take any $p < 1$ and color this random tessellation randomly, then the origin will a.s.\ be contained in a finite black cluster, as it is eventually enclosed by four white rectangles. Hence $p_c$ would be equal to one in this example.

\begin{figure}
  \centering
    \includegraphics[scale=.3]{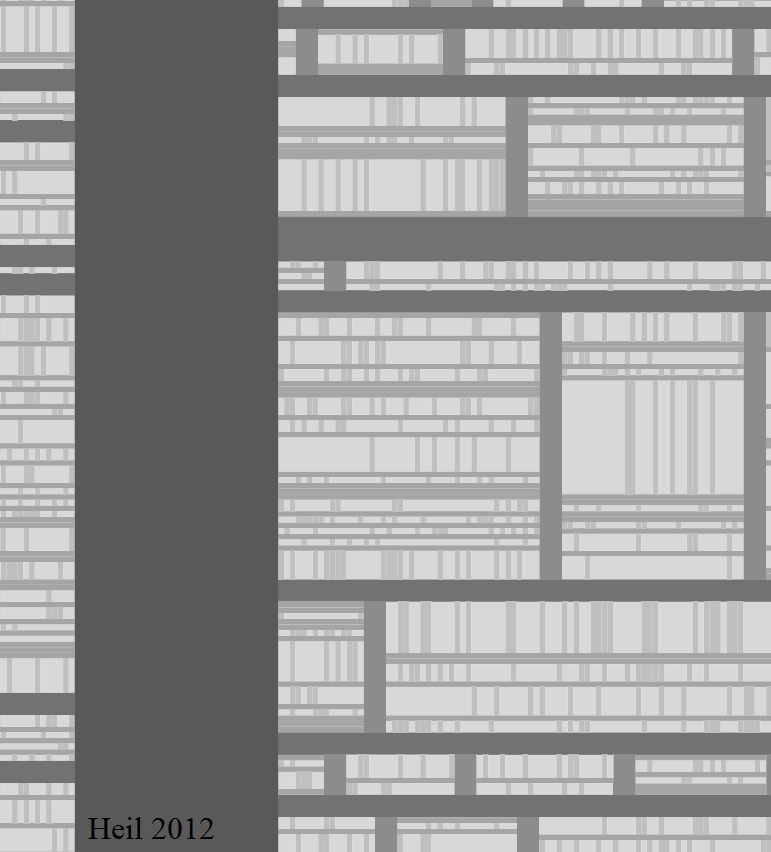}
  \caption{A sample of the streetgrid from \cite{heil2012stationary} where different cells are colored in different shades of grey.}
  \label{fig:streetgrid}
\end{figure}

\section{Framework I: Scale-mixing tessellations}\label{sec:scale-mixing_tess}
In the example above we have seen, that the standard mixing assumption does not imply a non-trivial phase transition. A natural question is, if there is a mixing condition, that is weaker than some continuous $k$-dependence analog and that ensures $p_c \in (0,1)$. The answer is yes and will be given in this section.

We say a function $f:\mathbf{M} \to \R$ is \emph{determined} by a set $A \subset \R^d$ if $f(m) = f(\tilde m)$ for any two $m, \tilde m \in \mathbf{M}$ with $m \cap A = \tilde m \cap A := \{Z \cap A \mid Z \in \tilde m, Z \cap A \neq \emptyset \}$. This means, if we know how $m$ looks in $A$, the value of $f(m)$ is fixed. A set $E \in \mathfrak{F}(\mathbf{M})$ is called determined by $A$, if $\ind_E$ is determined by $A$. We use this notion in the same way for point processes later on.

We remark, that a random tessellation is already mixing if \eqref{eq:mixing} holds for sets $A, A$ determined by cuboids $Q_1$ and $Q_2$ respectively \cite{schneider2008stochastic}. Hence heuristically we can say, that a random tessellation is mixing if events become more and more independent while we shift them away from each other. If we replace the shift in this heuristic by a scaling, we arrive at our new notion of scale-mixing.

A stationary random tessellation (or point process) $M$ is called \emph{scale-mixing} if for any two disjoint cuboids $Q := [a,b] := \times_{i = 1}^d [a_i, b_i]$ and $Q' := [a',b']$ with $a,a', b,b' \in \R^d$ we have
\begin{align}\label{eq:scale-mixing}
  \lim_{t\to \infty } \sup_{E_t, E'_t \in \mathfrak{F}(\mathbf{M})} |\P[M \in E_t \cap E'_t] - \P[M \in E_t] \P[M \in E'_t]| = 0
\end{align}
where $E_t$ and $E'_t$ are determined by $tQ$ and $tQ'$ respectively. For our applications, we have to ensure a certain speed of convergence in \eqref{eq:scale-mixing}. Therefor we say, that $M$ is \emph{scale-mixing of polynomial speed} (abbreviated by \emph{smp}) if for any two disjoint cuboids $Q, Q'$ there are constants $c_1, c_2 > 0$ such that for $t \in \R$ large enough and for all sets $E_t, E'_t \in \mathfrak{F}(\mathbf{M})$ determined by $t Q$ and $t Q'$ respectively, we have
\begin{align}\label{eq:def-smp}
  |\P[M \in E_t \cap E'_t] - \P[M \in E_t] \P[M \in E'_t]| \leq c_1 t^{-c_2}.
\end{align}

By replacing $M$ with $M_X$ and $\mathfrak{F}(\mathbf{M})$ with $\mathfrak{F}(\mathbf{M}_c)$, we obtain the corresponding definition of smp for independently marked stationary random tessellations. The first step now is to establish the link between smp for a random tessellation $M$ and its marked version $M_X$.

A function $z:\cC^d \to \R^d$ is called a center function if it is measurable and translation covariant, i.e.\ if $z(Z + x) = z(Z) + x$ for $x \in \R^d$, $Z \in \cC^d$. For the rest of the paper $z$ will be an arbitrary center function such that $z(Z) \in Z$. One might think of $z$ for example as the barycenter of $Z$ (further examples can be found in \cite{schneider2008stochastic}).

We recall the well known Campbell theorem (see \cite{schneider2008stochastic}), which says that for a stationary tessellation $M$ there is a measure $\Q$ concentrated on the cells $Z \in \cC^d$ with $z(Z) = 0$ and a number $\gamma \geq 0$ such that for any measurable $f:\R^d \times \cC^d \to [0,\infty )$ we have
\begin{align}\label{eq:campbell}
  \E\bigg[ \sum_{Z \in M} f(z(Z), Z-z(Z))\bigg] = \gamma \int_{\R^d} \int_{\cC^d} f(x,Z)\ \Q(dZ)\ dx.
\end{align}
The measure $\Q$ is called the \emph{distribution of the typical cell} and can be interpreted as the distribution of $Z - z(Z)$ if $Z$ is chosen ``uniformly from all cells of $M$''. The number $\gamma$ is the intensity of $M$ and is equal to the expected number of cell centers in the unit cube.

\begin{Lemma}
  Let $M$ be stationary tessellation, $c_1, c_2 > 0$ be functions of pairs of disjoint rectangles and $c_3, c'_3 > 0$. The stationary tessellation $M$ is smp with $c_1, c_2$ iff, for any two disjoint cuboids $Q$ and $Q'$ all large enough $t$ and all measurable $f:\mathbf{M} \to [0,c_3]$ determined by $tQ$ as well as $g:\mathbf{M} \to [0,c_3']$ determined by $tQ'$ we have
  \begin{align}\label{eq:smp-via-ev}
    |\E[f(M) g(M)] - \E[f(M)] \E[g(M)] | \leq c_3 c_3' c_1 t^{-c_2}.
  \end{align}
  Moreover, if the diameter of the typical cell has a finite $d + \varepsilon$ Moment, i.e.\ if there is an $\varepsilon > 0$ such that
  \begin{align}\label{eq:le41:moment-condition}
    \int_{\cC^d} \diam(Z)^{d+\varepsilon}\ \Q(dZ) < \infty,
  \end{align}
  then there is a constant $c_4 > 0$ such that for all $p \in [0,1]$ and $E_t, E'_t \in \mathfrak{F}(\mathbf{M}_c)$ determined by $tQ$ and $tQ'$ respectively, we have
  \[
    |\P_{p}[M_X \in E_t \cap E'_t] - \P_{p}[M_X \in E_t] \P_{p}[M_X \in E'_t] | \leq c_1 t^{-c_2} + c_4 t^{-\varepsilon}.
  \]
  Hence $M_X$ is smp.
\end{Lemma}
\emph{Proof:} The ``only if part'' is proved first for functions $f$ and $g$ of the form $\sum_{i = 1}^n h_i \ind_{H_i}$ with $h_i \in [0, 1/n]$ and $H_i$ determined by $tQ$ or $tQ'$ respectively. This is an easy exercise and standard approximation arguments finish this part. The other direction is trivial.

To prove the second assertion, we observe, that the map $f:m \mapsto \P_{p}[m_X \in E_t]$ is determined by $tQ$. By \eqref{eq:doubly-stochastic-model} we have $\P_{p}[M_X \in E_t] = \E[f(M)]$ and defining $g:m \mapsto \P_{p}[m_X \in E'_t]$ yields
\begin{align*}
  & \quad \ |\P_p[M_X \in E_t \cap E'_t] - \P_p[M_X \in E_t] \P_p[M_X \in E'_t]| \\
  & \leq |\P_p[M_X \in E_t \cap E'_t] - \E[f(M)g(M)]| + c_1 t^{-c_2}
\end{align*}
with the use of the first assertion and the triangle inequality.

Let $A_t$ be the set of tessellations that contain no cell that intersects $tQ$ as well as $tQ'$. It is clear, that if $m \in A_t$ then the coloring of $tQ$ is independent of the coloring of $tQ'$ and hence
\[
  \P_{p}[m_X \in E_t \cap E'_t] = \P_{p}[m_X \in E_t] \P_{p}[m_X \in E'_t].
\]
An easy calculation shows, that this implies
\[
  |\P_p[M_X \in E_t \cap E'_t] - \E[f(M)g(M)]| \leq \P[M \notin A_t].
\]
We will use the moment condition on the diameter of the typical cell to bound this probability.

Let $r \in \R$ be such that $Q, Q' \subset B_r$ and $\delta > 0$ such that $\delta r < \inf\{\|x-x'\|_2 \mid x \in Q, x' \in Q'\}$. If all cells of an $m \in \mathbf{M}$ that intersect $B_{tr}$ have a diameter less than $t \delta r$ then $m \in A_t$. Hence
\begin{align}
  \P[M \notin A_t] \leq \E\bigg[ \sum_{Z \in M} \ind\{Z \cap B_{tr} \neq \emptyset , \diam(Z) \geq t \delta r \}\bigg].
\end{align}
The use of Campbell's Formula \eqref{eq:campbell} together with basic estimates and the Markov-inequality yields
\begin{align*}
  \P[M \notin A_t] & \leq \gamma \int_{\R^d} \int_{\cC^d} \ind\{(Z+x) \cap B_{tr} \neq \emptyset , \diam(Z) \geq t \delta r\}\ \Q(dZ)\ dx \\
  & \leq \gamma \int_{B_{tr(1 + \delta)}} \int_{\cC^d} \ind\{\diam(Z) \geq t \delta r\} \ \Q(dZ)\ dx \\
  & \quad +  \gamma \int_{B_{tr(1 + \delta)}^c} \int_{\cC^d} \ind\{\diam(Z) \geq \|x\|_2 - t r\} \ \Q(dZ)\ dx \\
  & \leq \gamma \kappa_d (rt(1+\delta))^d t^{-d-\varepsilon } \int_{\cC^d} \diam(Z)^{d + \varepsilon}\ \Q(dZ) \\
  & \quad + \gamma d \kappa_d \int_{tr(1 + \delta)}^\infty y^{d-1} (y - tr)^{-d-\varepsilon} \int_{\cC^d} \diam(Z)^{d + \varepsilon}\ \Q(dZ)\ dy,
\end{align*}
where $\kappa_d$ is the volume of the unit ball in $\R^d$. Our moment condition \eqref{eq:le41:moment-condition} and an elementary evaluation of the remaining integral imply, that there is a $c_4 \in \R$ such that
\[
  \P[M \in A_t^c] \leq c_4 t^{-\varepsilon}
\]
which finishes the proof.\qed \bigskip

\subsection*{A non-trivial phase transition}
Now that we know some basic properties of smp tessellations the next theorem shows, that this property fits perfectly to show the existence of a non-trivial phase transition.

\begin{Thm}\label{thm:smp-implies-non-trivial-phasetransition}
  Let $M$ be an smp tessellation of $\R^d$, $d \geq 2$. If the diameter of the typical cell has a finite $d+\varepsilon$ moment for some $\varepsilon > 0$, then
  \[
    p_c(M) \in (0,1).
  \]
\end{Thm}
\emph{Proof:} We start by showing $p_c(M) < 1$. The idea for this part of the proof is to substitute the independence assumption in the second proof of Theorem 10 in \cite{bollobas2006percolation} by smp and relate the crossing probabilities of certain rectangles.

Let $M|_{\R^2} := \{Z \cap \R^2 \times \{0\}^{d-2} \mid Z \in M, Z \cap \R^2 \times \{0\}^{d-2} \neq \emptyset \}$ and observe that a.s.\ all cells in $M|_{\R^2}$ are 2-dimensional. Hence we can identify $M|_{\R^2}$ with a tessellation of $\R^2$ that is again smp and where the diameter of the typical cell fulfills at least the same moment conditions, as it only becomes smaller. If we find a $p < 1$ large enough, such that there is an infinite black cluster in $(M|_{\R^2})_X$, then there also exists an infinite black cluster in $M_X$. Hence we might restrict ourselves to the case $d = 2$.

For $a, b \in \R^d$ we will write again $[a,b]$ for the cuboid $\times_{i = 1}^d [a_i, b_i]$. For $a,b \in \R^2$ let $H(a,b)$ be the event, that there is a horizontal black crossing in $[a,b]$, i.e.\ there is curve in $[a,b]$ connecting $\{a_1\} \times [a_2,b_2]$ with $\{b_1\} \times [a_2,b_2]$ that uses only the interiors and common faces of black cells of $M_X$. In the same way we define vertical crossings $V(a,b)$ and remark, that the events $H(a,b)$ and $V(a,b)$ are determined by $[a,b]$.

\begin{figure}
  \centering
  \includegraphics{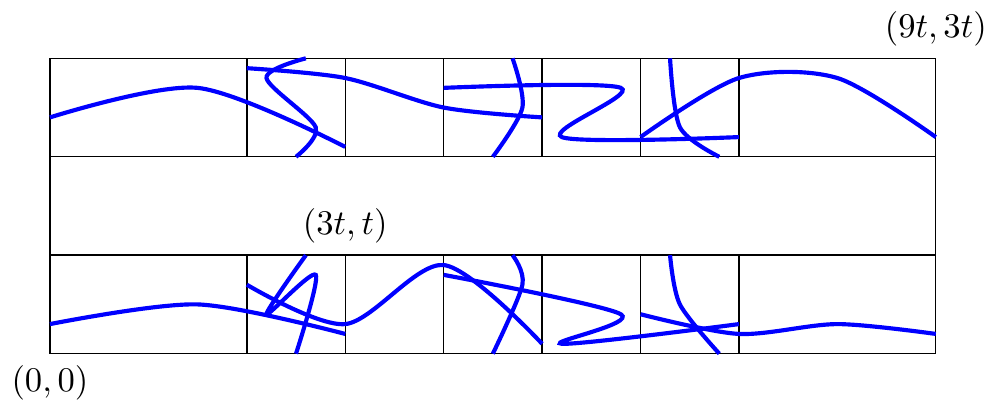}
  \caption{Horizontal and vertical crossings of certain small rectangles imply the crossing of the large rectangle (compare with \eqref{eq:skalierungsmischende_mosaike_nichttriviales_p_c:1} )}
  \label{fig:smp-p_c<1}
\end{figure}

\begin{figure}
  \centering
  \includegraphics{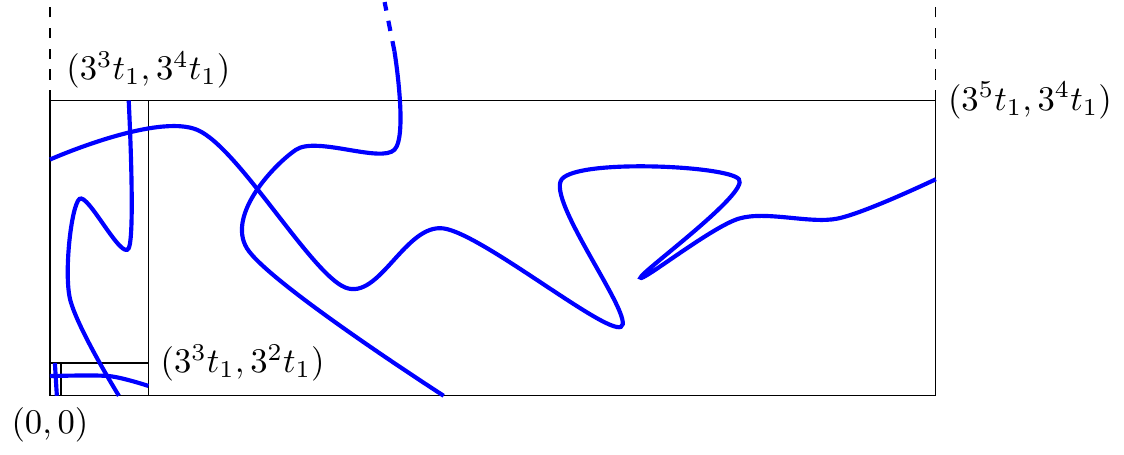}
  \caption{An infinite cluster evolves from rectangle crossings}
  \label{fig:infinite-cluster-from-crossings}
\end{figure}

We will now use the fact that $M_X$ is smp and stationary, to relate the probability of $H\left(\begin{psmallmatrix} 0 \\ 0 \\ \end{psmallmatrix}, \begin{psmallmatrix} 3t \\ t \\ \end{psmallmatrix}\right)$ for different values of $t$. Figure \ref{fig:smp-p_c<1} will make clear, that for $t$ large enough
\begin{align}
\begin{aligned}\label{eq:skalierungsmischende_mosaike_nichttriviales_p_c:1}
  &\quad \ \P_p\left[ H\left( \begin{psmallmatrix} 0 \\ 0 \\ \end{psmallmatrix}, \begin{psmallmatrix} 9t \\ 3t \\ \end{psmallmatrix}\right)^c \right] \\
  & \leq \P_p\left[ H\left( \begin{psmallmatrix} 0 \\ 0 \\ \end{psmallmatrix}, \begin{psmallmatrix} 9t \\ t \\ \end{psmallmatrix}\right)^c \cap  H\left( \begin{psmallmatrix} 0 \\ 2t \\ \end{psmallmatrix}, \begin{psmallmatrix} 9t \\ 3t \\ \end{psmallmatrix}\right)^c \right] \\
  & \leq \P_p\left[ H\left( \begin{psmallmatrix} 0 \\ 0 \\ \end{psmallmatrix}, \begin{psmallmatrix} 9t \\ t \\ \end{psmallmatrix}\right)^c \right]^2 + c_1 t^{-c_2} \\
  & \leq \P_{p}\left[ H\left( \begin{psmallmatrix} 0 \\ 0 \\ \end{psmallmatrix}, \begin{psmallmatrix} 3t \\ t \\ \end{psmallmatrix}\right)^c
  \cup V\left( \begin{psmallmatrix} 2t \\ 0 \\ \end{psmallmatrix}, \begin{psmallmatrix} 3t \\ t \\ \end{psmallmatrix} \right)^c
  \cup H\left( \begin{psmallmatrix} 2t \\ 0 \\ \end{psmallmatrix}, \begin{psmallmatrix} 5t \\ t \\ \end{psmallmatrix}\right)^c
  \cup V\left( \begin{psmallmatrix} 4t \\ 0 \\ \end{psmallmatrix}, \begin{psmallmatrix} 5t \\ t \\ \end{psmallmatrix} \right)^c \right.\\
  & \hspace{2cm} \left.
  \cup H\left( \begin{psmallmatrix} 4t \\ 0 \\ \end{psmallmatrix}, \begin{psmallmatrix} 7t \\ t \\ \end{psmallmatrix}\right)^c
  \cup V\left( \begin{psmallmatrix} 6t \\ 0 \\ \end{psmallmatrix}, \begin{psmallmatrix} 7t \\ t \\ \end{psmallmatrix} \right)^c
  \cup H\left( \begin{psmallmatrix} 6t \\ 0 \\ \end{psmallmatrix}, \begin{psmallmatrix} 9t \\ t \\ \end{psmallmatrix}\right)^c
  \right]^2 + c_1 t^{-c_2} \\
  & \leq \big( 4 \P_p\left[ H\left( \begin{psmallmatrix} 0 \\ 0 \\ \end{psmallmatrix}, \begin{psmallmatrix} 3t \\ t \\ \end{psmallmatrix}\right)^c \right] + 3 \P_p\left[ V\left( \begin{psmallmatrix} 0 \\ 0 \\ \end{psmallmatrix}, \begin{psmallmatrix} t \\ t \\ \end{psmallmatrix} \right)^c \right]\big)^2 + c_1 t^{-c_2} \\
  & \leq \big( 7 \max\big\{ \P_p\left[ H\left( \begin{psmallmatrix} 0 \\ 0 \\ \end{psmallmatrix}, \begin{psmallmatrix} 3t \\ t \\ \end{psmallmatrix}\right)^c \right], \P_p\left[ V\left( \begin{psmallmatrix} 0 \\ 0 \\ \end{psmallmatrix}, \begin{psmallmatrix} t \\ 3t \\ \end{psmallmatrix} \right)^c \right]  \big\}\big)^2 + c_1 t^{-c_2},
\end{aligned}
\end{align}
where $c_1, c_2$ are the smp constants corresponding to the rectangles $\left[ \begin{psmallmatrix} 0 \\ 0 \\ \end{psmallmatrix}, \begin{psmallmatrix} 9 \\ 1 \\ \end{psmallmatrix} \right]$ and $\left[ \begin{psmallmatrix} 0 \\ 2 \\ \end{psmallmatrix}, \begin{psmallmatrix} 9 \\ 3 \\ \end{psmallmatrix} \right]$. The same relation holds for $V\left( \begin{psmallmatrix} 0 \\ 0 \\ \end{psmallmatrix}, \begin{psmallmatrix} 3t \\ 9t \\ \end{psmallmatrix}\right)$ and we obtain constants $c_3, c_4 > 0$ and a $t_0$ that depend on the smp constants of $M$ such that
\begin{align}
  f_{n+1} \leq 49 f_n^2 + c_3 (3^n t_1)^{-c_4}, \quad \forall n \in \N,\ t_1 \geq t_0
\end{align}
where
\[
  f_n := \max \left\{ \P_p\left[ H\left( \begin{psmallmatrix} 0 \\ 0 \\ \end{psmallmatrix}, \begin{psmallmatrix} 3^{n+1} t_1 \\ 3^n t_1 \\ \end{psmallmatrix}\right)^c \right], \P_p\left[ V\left( \begin{psmallmatrix} 0 \\ 0 \\ \end{psmallmatrix}, \begin{psmallmatrix} 3^n t_1 \\ 3^{n+1} t_1 \end{psmallmatrix} \right)^c \right]  \right\}.
\]
A simple induction shows, that if we choose $t_1$ large enough such that $c_3 t_1^{-c_4} < (4 \cdot 49 \cdot 3^{c_4})^{-1}$ and if $f_1 \leq (2 \cdot 49 \cdot 3^{c_4})^{-1}$ then $f_n \leq 3^{-c_4 n}$. The value of $f_1$ can be made arbitrary small, by choosing $p$ large enough, as there certainly is a crossing in $[a,b]$ if all cells that intersect $[a,b]$ are black. Hence the Borel-Cantelli lemma implies that a.s.\ only a finite number of the events
\[
   H\left( \begin{psmallmatrix} 0 \\ 0 \\ \end{psmallmatrix}, \begin{psmallmatrix} 3^{2n+1} t_1 \\ 3^{2n} t_1 \\ \end{psmallmatrix}\right),\ V\left( \begin{psmallmatrix} 0 \\ 0 \\ \end{psmallmatrix}, \begin{psmallmatrix} 3^{2n+1} t_1 \\ 3^{2n+2} t_1 \\ \end{psmallmatrix}\right), \quad  n \in \N
\]
won't hold. This yields the existence of an infinite black path in $\cG_{M, X}$, as the horizontal crossing in $\left[  \begin{psmallmatrix} 0 \\ 0 \\ \end{psmallmatrix}, \begin{psmallmatrix} 3^{2n+1} t_1 \\ 3^{2n} t_1 \\ \end{psmallmatrix} \right]$ and the vertical crossing in  $\left[  \begin{psmallmatrix} 0 \\ 0 \\ \end{psmallmatrix}, \begin{psmallmatrix} 3^{2n+2} t_1 \\ 3^{2n+1} t_1 \\ \end{psmallmatrix} \right]$ have to intersect (see Figure \ref{fig:infinite-cluster-from-crossings}).

The argument for $p_c(M) > 0$ is similar, but the geometry is a bit more complicated. We consider the cuboids $Q_i := [0,3]^{i-1} \times [0,1] \times [0,3]^{d-i}$, $i \in [d]$ and the sets of colored tessellations $A(i,t) \in \mathfrak{F}(\mathbf{M}_c)$ where $tQ_i$ is crossed in the short direction, i.e.\ there is a black curve in $tQ_i$ starting in $[0,3t]^{i-1} \times \{0\} \times [0,3t]^{d - i}$ and ending in $[0,3t]^{i-1} \times \{t\} \times [0,3t]^{d - i}$.

If a curve starts in $t[-\tfrac{1}{2}, \tfrac{1}{2}]^d$ and ends in $t[-\tfrac{3}{2}, \tfrac{3}{2}]^d$ then at least one of the $2d$ cuboids $tQ_i + \tfrac t2 \mathbf{e}_i + \sum_{j \in [d]\setminus \{i\}} \tfrac{3t}{2} \mathbf{e}_j$, $tQ_i - \sum_{j \in [d]} \tfrac{3t}{2} \mathbf{e}_j$, $i \in [d]$ that are arranged ``around'' $t[-\tfrac{1}{2}, \tfrac{1}{2}]^d$, will be crossed in the short direction, by this curve.

\begin{figure}
  \centering
  \includegraphics{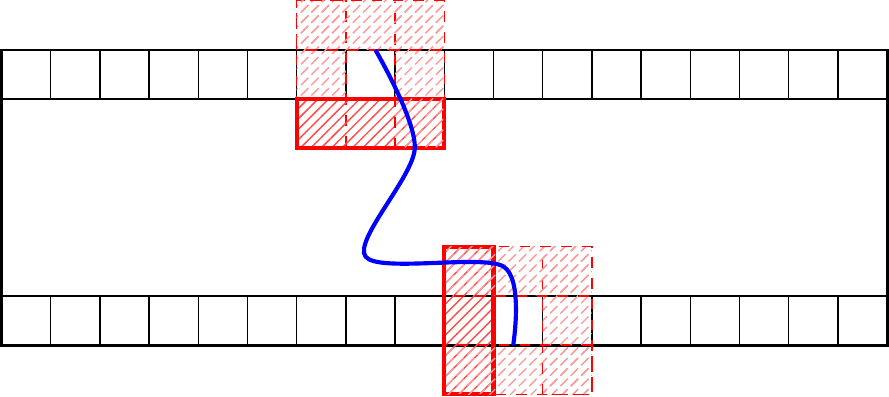}
  \caption{The blue crossing of the large rectangle implies that at least one square on the upper side and one square on the lower side are ``left'' by this path. Hence there are at least two disjoint red rectangles that are crossed in the short direction.}
  \label{fig:smp-p_c>0}
\end{figure}

We want to construct a recursion for the probability that $M_X \in A(i,t)$. To this end we consider $6t Q_i$ and cover both of its sides $[0,3 \cdot 6t]^{i-1} \times \{0\} \times [0,3 \cdot 6 t]^{d - i}$ and $[0,3 \cdot 6 t]^{i-1} \times \{t\} \times [0,3 \cdot 6 t]^{d - i}$ with a finite number of cubes of sidelength $t$ with non-overlapping interiors, that are all contained in $6t Q_i$. If $M_X \in A(i,6t)$ holds, then on each of the two sides there has to be a cube in which a curve starts that ends further than $t$ away from it. This implies, that there is a translation of one of the cuboids $tQ_i$ that is crossed in the short direction (see Figure \ref{fig:smp-p_c>0}).

Now we may do a similar calculation as in the first part of the proof to obtain constants $c_1, c_2, c_3, t_0 \in (0,\infty )$ such that for all $t_1 \geq t_0$
\[
  \max_{i \in [d]}\P_{p}[M_X \in A(i,t_1 )] \leq c_3 \max_{i \in [d]}\P_{p}[M_X \in A(i, t_1)]^2 + c_1 (6^n t_1)^{-c_2}.
\]

A cuboid won't be crossed by a black curve if all cells that intersect this cuboid are white. Hence we find a $t_1 \in \R$ and a $p > 0$ such that
\[
  \max_{i \in [d]} \P_{p}[M_X \in A(i,6^n t_1)] \leq 6^{-c_2 n}, \quad n \in \N.
\]
The Borel-Cantelli lemma ensures now that the zero cell is contained in a finite black cluster a.s.\ .\qed\bigskip

The proof may also be a blueprint for other percolation models in euclidian space. One only has to verify the smp condition and make sure that crossings of cuboids are very probable or improbable. We also want to remark that it is easy to show, that the clustervolume has a polynomial tail, if $\max_{i \in [d]} \P_{p}[M_X \in A(i,6^n t_1)] \leq 6^{-c_2 n}$, $n \in \N$.

It can be observed in the proof, that we only need a tiny part of the smp property namely that it holds for crossing events. We could imagine that there are models, where one could verify just this part, to proof a non-trivial phase transition.

\subsection*{Examples}
The bad news first: STIT tessellations and Poisson hyperplane tessellations are not scale-mixing (For an introduction to STIT we refer to \cite{nagel2005crack}, material on Poisson hyperplane tessellations can be found in \cite{schneider2008stochastic}). We only sketch the arguments but won't go into the details.

To simplify the argument for STIT, we might think of a 2-dimensional version with only horizontal and vertical lines, that have equal probability in the directional distribution. Now we fix two unit squares $Q_1, Q_2 := Q_1 + \binom{1.5}{0}$ and consider the two events ``$tQ_i$ has a horizontal line from left to right''. It is easy to show that the difference between the probability of the intersection of these two events and the product of their probabilities is bounded from below by some constant $c > 0$. Most likely other counterexamples can be found for other directional distributions or in higher dimensions, but we won't elaborate on that. A similar counterexample works for any directional distribution in the Poisson hyperplane tessellation.

Fortunately there are some interesting results for Voronoi tessellations induced by point processes. Let $\varphi \subset \R^d$ be a locally finite set of points. The set
\begin{align}
  Z(x,\varphi) := \{y \in \R^d \mid \forall z \in \varphi:\ \|y - x\|_2 \leq \|y - z\|_2 \},\quad x \in \varphi
\end{align}
is called the \emph{Voronoi cell} of $x$ in $\varphi$ and the collection
\begin{align}
  \cV(\varphi) := \{Z(x,\varphi) \mid x \in \varphi\}
\end{align}
is called the \emph{Voronoi tessellation} induced by $\varphi$. If we replace $\varphi$ by a stationary point process $\Phi$ on $\R^d$ we obtain the random tessellation $\cV(\Phi)$. This and more related results can be found in \cite{schneider2008stochastic}.

In the remainder of this section, we will investigate which classes of point processes induce smp Voronoi tessellations. The following preparing lemma will show that the void probabilities of smp point processes decay at least polynomial.

\begin{Lemma}\label{le:void-probabilities-smp}
  Let $\Phi$ be a smp point process on $\R^d$ with intensity $\gamma > 0$. If $Q$ is a cuboid, then there are constants $c_1, c_2 > 0$ such that
  \begin{align}
    \P[\Phi(tQ) = 0] \leq c_1 t^{-c_2}
  \end{align}
  for large enough $t$
\end{Lemma}
\emph{Proof:} Similar to the proof of Theorem \ref{thm:smp-implies-non-trivial-phasetransition} we derive a relation between the events ``$\Phi(3tQ) = 0$'' and ``$\Phi(tQ) = 0$''. It is clear, that we can find $x_1, x_2 \in \R^d$ such that the cuboids $tQ + x_1$ and $tQ + x_2$ are disjoint and are contained in $3tQ$. Hence there are constants $c_3,c_4 > 0$ such that for $t$ large enough
\[
  \P[\Phi(3tQ) = 0] \leq \P[\Phi(tQ + x_1) = 0, \Phi(tQ + x_2) = 0] \leq \P[\Phi(tQ) = 0]^2 + c_3 t^{-c_4}.
\]
As $\gamma > 0$ we have that
\[
  \P[\Phi(tQ) = 0] \tofor{t} 0
\]
and an easy induction shows the existence of a $t_1 \in \R$ and a constant $c_5 > 0$ such that
\[
  \P[\Phi(3^n t_1 Q) = 0] \leq c_5 3^{-c_4 n}, \quad n \in \N.
\]
Monotonicity implies, that $\P[\Phi(3^x 3 t_1 Q) = 0] \leq c_5 3^{-c_4 x}$ for all $x > 1$.\qed\bigskip

\begin{figure}
  \centering
  \includegraphics{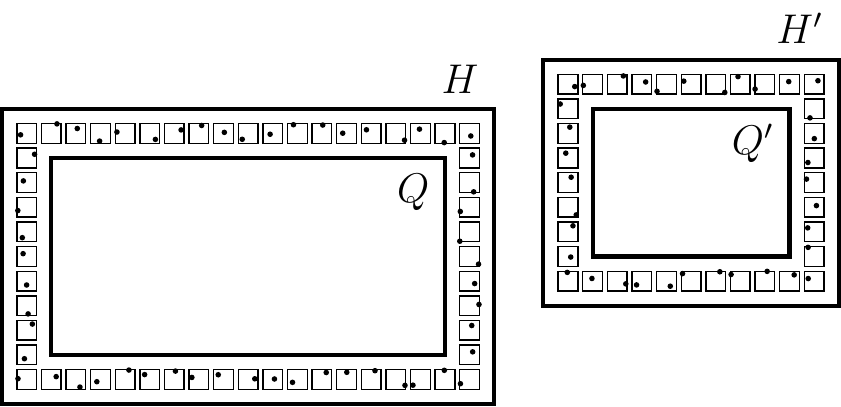}
  \caption{Example for a possible choice of all the cuboids.}
  \label{fig:smp-pp-implies-smp-voronoi}
\end{figure}

\begin{Thm}\label{thm:smp-pp-implies-smp-voronoi}
  If $\Phi$ is an smp point process on $\R^d$, then $\cV(\Phi)$ is smp too.
\end{Thm}
\emph{Proof:} Let $Q$ and $Q'$ be two disjoint cuboids. We choose another two disjoint cuboids $H$ and $H'$ such that $Q$ and $Q'$ lie in the interior of $H$ and $H'$ respectively. We argue that we can choose some more cuboids $W_1,\dots, W_n \subset H \setminus Q$ and $W'_1,\dots,W'_{n'} \subset H' \setminus Q'$ which are translates of one another and have the following property. If each of these cuboids contains one point of $\Phi$, then $\cV(\Phi)$ restricted to $Q$ ($Q'$) is determined by the points in $H$ ($H'$). This is due to the fact that these cuboids can be chosen so small and close to $Q$ that no point in $Q$ has it's closest point in $\Phi$ outside of $H$ (see Figure \ref{fig:smp-pp-implies-smp-voronoi}). This doesn't change, if we scale everything by a factor $t$.

Let $E \in \cN(\R^d)$ be the set of point configurations, where each of the cubes $tW_i$, $i \in [n]$ contains at least one point. Let $E' \in \cN(\R^d)$ be the same set for $tW'_j$, $j \in [n']$. By Lemma \ref{le:void-probabilities-smp} and subadditivity we have that
\begin{align}\label{eq:smp-pp-implies-smp-voronoi:-1}
  \P[\Phi \in E^c] \leq n c_1 t^{-c_2}
\end{align}
and
\begin{align}\label{eq:smp-pp-implies-smp-voronoi:0}
  \P[\Phi \in E'^c] \leq n' c_1 t^{-c_2}
\end{align}
for $t$ large enough and suitable constants $c_1, c_2 > 0$.

Now we check \eqref{eq:def-smp}. By the triangle inequality and trivial estimates we have for $A, A' \in \mathfrak{F}(\mathbf{M}_c)$ determined by $tQ$ and $tQ'$ respectively as well as large enough $t$ that
\begin{align}
\begin{aligned}\label{eq:smp-pp-implies-smp-voronoi:1}
  &\quad \ | \P_p[\cV(\Phi)_X \in A \cap A'] - \P_p[\cV(\Phi)_X \in A] \P_p[ \cV(\Phi)_X \in A'] |\\
  &\leq | \P_p[\cV(\Phi)_X \in A, \Phi \in E,\ \cV(\Phi)_X \in A', \Phi \in E' ] \\
  & \hspace{1cm} - \P_p[\cV(\Phi)_X \in A, \Phi \in E ] \P_p[ \cV(\Phi)_X \in A', \Phi \in E' ] |\\
  & \hspace{1cm} + 2 \P[E^c] + 3 \P[E'^c] .
\end{aligned}
\end{align}
The construction of $E$ ensures, that no cell that intersects $tQ$ will intersect $tQ'$ and vice versa. Hence by \eqref{eq:doubly-stochastic-model}
\begin{align}
\begin{aligned}\label{eq:smp-pp-implies-smp-voronoi:2}
  &\quad \ \P_p[\cV(\Phi)_X \in A, \Phi \in E,\ \cV(\Phi)_X \in A', \Phi \in E' ] \\
  & = \int \P_{p}[\cV(\varphi)_X \in A, \varphi \in E,\ \cV(\varphi)_X \in A', \varphi \in E' ] \ \P_{\Phi}(d\varphi) \\
  & = \int \P_{p}[\cV(\varphi)_X \in A, \varphi \in E] \P_{p}[\cV(\varphi)_X \in A', \varphi \in E' ] \ \P_{\Phi}(d\varphi) .
\end{aligned}
\end{align}

The construction of $E$ also ensures that if a function $f:\mathbf{M} \to [0,1]$ is determined by $tQ$, then, the function $g:\mathbf{N} \to [0,1],\ \varphi\mapsto f(\varphi) \ind_E(\varphi)$ is determined by $tH$. Hence we may apply \eqref{eq:smp-via-ev} to the function $\varphi\mapsto \P_{p}[\cV(\varphi)_X \in A, \varphi \in E]$ and its prime pendant which are determined by $tH$ and $tH'$ respectively. With the help of \eqref{eq:smp-pp-implies-smp-voronoi:2} we obtain that
\begin{align}
\begin{aligned}\label{eq:smp-pp-implies-smp-voronoi:3}
  &\quad \ | \P_p[\cV(\Phi)_X \in A, \Phi \in E,\ \cV(\Phi)_X \in A', \Phi \in E' ] \\
  & \hspace{1cm} - \P_p[\cV(\Phi)_X \in A, \Phi \in E ] \P_p[ \cV(\Phi)_X \in A', \Phi \in E' ] |\\
  & \leq c_3 t^{-c_4},
\end{aligned}
\end{align}
where $c_3, c_4 > 0$ are again suitable constants following from the fact that $\Phi$ is smp. Combining \eqref{eq:smp-pp-implies-smp-voronoi:-1}, \eqref{eq:smp-pp-implies-smp-voronoi:0}, \eqref{eq:smp-pp-implies-smp-voronoi:1} and \eqref{eq:smp-pp-implies-smp-voronoi:3} proves the  assertion.\qed \bigskip

Theorem \ref{thm:smp-pp-implies-smp-voronoi} transfers the problem of finding smp tessellations to finding smp point processes. It is clear, that any $k$-dependent point process is smp which includes the Poisson process and some perturbed lattice point processes. Two more interesting classes of point processes are Poisson cluster and Gibbsian point processes. We will show, that large subclasses of them are smp.

We start with the Poisson cluster process. Let $\Phi_0 = \sum_{n \in \N} \delta_{\zeta_n}$ be a Poisson process with intensity $\gamma_0 > 0$ on $\R^d$ and let $(\Psi_i)_{i \in \N}$ be an i.i.d. sequence of point processes that are independent of $\Phi_0$, have the common finite intensity measure $\Gamma_1$ and the distribution $\P_\Psi$. We call the point process
\[
  \Phi := \sum_{n \in \N} \zeta_n + \Psi_n
\]
the \emph{Poisson cluster process} with ground intensity $\gamma_0$ and cluster distribution $\P_\Psi$. This is a stationary point process with Intensity $\gamma_0 \Gamma_1(\R^d)$. The idea behind this construction is, that each point $\zeta_n$ of the Poisson process is replaced by a whole cluster of points $\Psi_n + \zeta_n$.

\begin{Thm}\label{thm:poisson-cluster-is-smp}
  Let $\Phi$ be a Poisson cluster process. If there are constants $c_1, c_2 > 0$ such that
  \[
    \P[\diam(\Psi_1) \geq r] \leq c_1 r^{-(d + c_2)}, \quad r > 0
  \]
  then $\Phi$ is smp.
\end{Thm}
\emph{Proof:} If each point of a Poisson process is shifted independently according to some fixed law, then the result is again a Poisson process of the same intensity hence we may assume, that $\P_{\Psi}$ is such that a.s. $\Psi_1 = 0$ or $\Psi_1(\{0\}) \geq 1$, i.e.\ either a cluster is empty or it has a point in the origin.

Let $Q$ and $Q'$ be two disjoint cuboids. We choose another two disjoint cuboids $H, H'$ and $c_3 > 0$ such that $Q + B_{c_3} \subset H$ and $Q' + B_{c_3} \subset H'$. For a Borel set $D$ we define
\begin{align}
  \Phi_D := \sum_{n \in \N} \ind\{\zeta_n \in D\} (\Psi_n + \zeta_n)
\end{align}
the cluster process that stems from $\Phi_0$ restricted to $D$. Related to this we define for disjoint Borel sets $D_1$, $D_2$ the event $E(D_1, D_2)$ which holds if no cluster $\Psi_n + \zeta_n$ with $\zeta_n \in D_1$ intersects $D_2$. It is easy to see, that $E(D_1 \cup \tilde D_1, D_2) = E(D_1, D_2) \cap E(\tilde D_1,D_2)$ if $\tilde D_1$ doesn't intersect $D_2$ either.

The definition of $E(D_1, D_2)$ implies, that if $A \in \cN(\R^d)$ is determined by $Q$ then
\begin{align}
  \{\Phi \in A\} \cap E(H^c, Q) = \{\Phi_{H} \in A\} \cap E(H^c, Q).
\end{align}
It is also clear, that the event $E(D_1, D_2)$ depends only on the Poisson points $\zeta_n \in D_1$ and the corresponding clusters. Hence $E(D_1, D_2)$, $E(D_3,D_4)$ and $\{\Phi_{D_5} \in A\}$ are independent events as long as $D_1$, $D_3$ and $D_5$ are disjoint. Taking all these arguments into account yields that
\begin{align}
\begin{aligned}\label{eq:poisson-cluster-smp:3}
  &\quad \ \P[ \{\Phi \in A \cap A'\} \cap E(H^c, Q) \cap E(H'^c, Q') ] \\
  & = \P[\{\Phi_H \in A\} \cap  \{\Phi_H' \in A'\} \cap E(H^c, Q) \cap E(H'^c, Q')] \\
  & = \P[\{\Phi_H \in A\} \cap E(H, Q') \cap \{\Phi_H' \in A'\} \cap E(H', Q) \\
  & \hspace{1cm} \cap E( (H \cup H')^c , Q) \cap E( (H \cup H')^c , Q')] \\
  & = \P[\{\Phi_H \in A\} \cap E(H, Q')] \P[\{\Phi_H' \in A'\} \cap E(H', Q)] \\
  & \hspace{1cm} \P[E( (H \cup H')^c , Q \cup Q')]
\end{aligned}
\end{align}
if $A' \in \cN(\R^d)$ is determined by $Q'$. Replacing $Q$, $Q'$, $H$ and $H'$ by their scaled version and applying similar arguments as in Theorem \ref{thm:smp-pp-implies-smp-voronoi}, we deduce from \eqref{eq:poisson-cluster-smp:3} that
\begin{align}
\begin{aligned}\label{eq:poisson-cluster-smp:4}
  & \quad \ |\P[\Phi \in A \cap A'] - \P[\Phi \in A] \P[\Phi \in A']| \\
  & \leq c_3 \P[E(tH, tQ')^c] + c_4 \P[E(tH',tQ )^c] + c_5 \P[E( (tH \cup tH')^c , tQ \cup tQ')^c]
\end{aligned}
\end{align}
for suitable $c_3,c_4,c_5 \in \N$.

Hence our final task is to show, that $\P[E(tH^c, tQ)^c]$ decays at least polynomial in $t$ as it dominates all three error terms in \eqref{eq:poisson-cluster-smp:4}. We have
\begin{align*}
  \P[E(tH^c, tQ)^c] & \leq \P[\exists n \in \N:\ \zeta_n \in tH^c, \diam(\Psi_n) \geq d(\zeta_n, tQ)] \\
  & \leq \E\left[ \sum_{n \in \N} \ind\{\zeta_n \in tH^c, \diam(\Psi_n) \geq d(\zeta_n, tQ)\}\right]
\end{align*}
where $d(\zeta_n,tQ) := \min_{y \in tQ} \|x - y\|_2$. The Mecke equation implies that
\[
  \P[E(tH^c, tQ)^c] \leq \gamma_0 \int_{tH^c} \P[\diam(\Psi_0) \geq d(x, tQ)]\ dx.
\]
Applying our assumption on the tail behaviour of the diameter of the clusters we get
\begin{align*}
  \P[E(tH^c, tQ)^c] & \leq \gamma_0 c_1 \int_{tH^c} d(x, tQ)^{-d-c_2}\ dx. \\
  & = \gamma_0 c_1 t^d \int_{H^c} d(ty, tQ)^{-d-c_2}\ dy \\
  & = \gamma_0 c_1 t^{-c_2} \int_{H^c} d(y, Q)^{-d-c_2}\ dy.
\end{align*}
An easy calculation shows, that the integral is finite which completes the proof.\qed\bigskip

The second class of smp point processes are the Gibbsian point processes studied in \cite{schreiber2013}. The underlying perfect simulation procedure ensures, that events determined by disjoined regions are independent as long as their so-called ancestor clans don't intersect. The diameter of these ancestor clans has an exponential tail and calculations like the one above result in an estimate like
\[
  |\P[\Phi \in A \cap A'] - \P[\Phi \in A] \P[\Phi \in A']| \leq c_1 e^{-t c_2}
\]
with a Gibbs process $\Phi$ from one of the classes defined in the introduction of \cite{schreiber2013}, suitable constants $c_1, c_2 > 0$, large enough $t$ and sets $A, A'$ determined by cuboids $tQ$ and $tQ'$ respectively (see Lemmas 3.4 and 3.5 in \cite{schreiber2013}).

\section{Framework II: Tame tessellations}\label{sec:tame_tess}
We have seen in the example before Section \ref{sec:scale-mixing_tess}, that very large cells might be a problem for a non-trivial phase transition and one could easily imagine, that too many very small cells might not be nice either. So, our second approach to get sufficient conditions for a non-trivial phase transition is, to bound the occurrence of too large or too small cells. We will need the notion of (greedy) lattice animals for this.

Let $\cG = (V, E)$ be a graph with root $\mathbf{0}$. The set $\cA(\cG)$ of \emph{graph animals} of $\cG$ contains all subsets of $V$ that contain the root and are connected in $\cG$. We define the set of animals of a given size $\cA_n(\cG) := \{\alpha \in \cA(\cG) \mid |\alpha| = n\}$, $n \in \N$. In the special case, where $\cG = \cZ^d := (\Z^d, \{\{u,v\} \subset \Z^d \mid \|u - v\|_2 = 1\})$ we write $\cA^{(d)} := \cA(\cZ^d)$ and call the elements \emph{lattice animals}.

To quantify how many large and small cells are too many, we use the notion of greedy lattice animals, that was introduced in \cite{cox1993greedy} and \cite{gandolfi1994greedy}, for two auxiliary random fields, that are constructed deterministically from our random tessellation. Cox et al.\ investigated the behaviour of
\[
  \max_{\alpha \in \cA^{(d)}} \frac{1}{n} \sum_{v \in \alpha} Y_v, \quad n \in \N
\]
(the greediest animal) where $Y = (Y_v)_{v \in \Z^d}$ is an i.i.d. random field. Among other results they found a condition on the tail behaviour of $Y_0$ such that the $\limsup_{n\to\infty}$ of the above maximum is almost surely finite.

\begin{figure}
  \centering
  \includegraphics{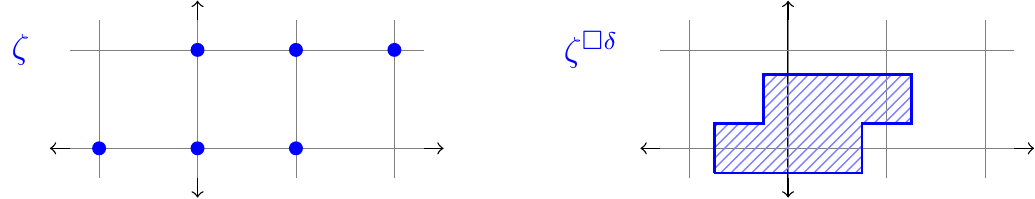}
  \caption{A set $\zeta \subset \Z^2$ and the corresponding $\zeta^{\square \delta}$ with $\delta = \tfrac{1}{2}$.}
  \label{fig:square_delta}
\end{figure}

We want to work with grids of various sizes and the corresponding cubes in this section. Therefor we denote by
\[
  \zeta^{\square \delta} := \delta \left( \zeta + [-\tfrac 12, \tfrac 12]^d \right)
\]
the union of boxes around a point or a set of points $\zeta \subset \Z^d$ in a grid with width $\delta$ (see Figure \ref{fig:square_delta}).

For a stationary random tessellation $M$ and $\delta > 0$ we define the auxiliary random fields $Y := (Y_v)_{v \in \Z^d} := (Y_v(M, \delta))_{v \in \Z^d}$ and $U := (U_v)_{v \in \Z^d} := (U_v(M, \delta))_{v \in \Z^d}$ by
\begin{align}
  \begin{aligned}
    Y_v & := |\{Z \in M \mid z(Z) \in v^{\square \delta}\}| \\
    U_v & := \ind\{\text{a cell of $M$ intersects $v^{\square \delta}$ and $\{w \in \Z^d \mid \|w - v\|_\infty \geq 2\}^{\square \delta}$}\}.
  \end{aligned}
\end{align}
As $M$ is stationary, $Y$ and $U$ are also stationary and there is a.s.\ no cell center on the boundary of a box. Hence a.s.\ no centers are counted twice by $Y$.

\begin{Def}
  An ergodic random tessellation $M$ is called \emph{tame} if
  \begin{enumerate}
    \item[(T1)] there is a $\delta_1 > 0$ and a $c_1 \in \R$ such that
        \[
          \limsup_{n\to \infty } \max_{\alpha \in \cA_n^{(d)}} \frac{1}{n} \sum_{v \in \alpha} Y_v(M, \delta_1) \leq c_1,
        \]
    \item[(T2)] there is a $\delta_2 > 0$ and a $c_2 < 1$ such that
        \[
          \limsup_{n\to \infty } \max_{\alpha \in \cA_n^{(d)}} \frac{1}{n} \sum_{v \in \alpha} U_v(M, \delta_2) \leq c_2.
        \]
  \end{enumerate}
\end{Def}

Occasionally it is the case, that $Y_0$ is a.s.\ bounded for some $\delta_1 > 0$ (e.g. if $M$ is the Voronoi-tessellation of some hard-core point process). Then (T1) holds trivially.

We want to avoid the use of two different grid sizes $\delta_1$ and $\delta_2$. Any cube $v^{\square \delta_1}$ is covered by at most $c_3 := c_3(\delta_1, \delta_2,d) \in \N$ cubes $w^{\square \delta_2}$. Hence for any $\alpha \in \cA_n^{(d)}$ there is a $\tilde \alpha \in \cA_{c_3 n}^{(d)}$ such that $\sum_{v \in \alpha} Y_v(M, \delta_2) \leq \sum_{v \in \tilde \alpha} Y_v(M, \delta_1)$. Hence if (T1) holds for some $\delta_1 > 0$ it also holds if we chose $\delta_1 := \delta_2$. Therefor we will always use $\delta$ instead of $\delta_1$ or $\delta_2$ from now on. We will write $\zeta^{\square} := \zeta^{\square \delta}$.

The tameness of a random tessellation ensures a non-trivial phase transition, as the next theorem will show, but it is also very useful when considering first passage percolation on random tessellations (see \cite{ziesche2016first}).

\begin{Thm}\label{thm:non-trivial-phase-transition}
  Let $M$ be a random tessellation of $\R^d$. If $M$ is tame, then there is a non-trivial phase transition. Moreover there are constants $c_1, c_2 \in \R$ such that a.s.\ for $n \in \N$ large enough
  \begin{align}\label{eq:non-trivial-phase-transition:claim1}
    |\cA_n(\cG_M)| \leq c_1^n
  \end{align}
   and
  \begin{align}\label{eq:non-trivial-phase-transition:claim2}
    |\partial^+ \alpha| \leq c_2 n,\quad  \alpha \in \cA_n(\cG_M).
  \end{align}
\end{Thm}
\emph{Proof:} We start with the proof of the animal bound, as it immediately implies $p_c \geq 1/c_1$ due to Lemma \ref{le:animal-bound-implies-p_c>0}.

The random tessellation $M$ is tame and hence there is a $\delta > 0$, $c_3 < 1$ and $c_4 \in \R$ such that a.s.\ there is an $n_0 \in \N$ with
\begin{align}\label{eq:non-trivial-phase-transition:1}
  \max_{\zeta \in \cA_n^{(d)}} \sum_{v \in \zeta} U_v \leq c_3 n, \quad n \geq n_0
\end{align}
and
\begin{align}\label{eq:non-trivial-phase-transition:2}
  \max_{\zeta \in \cA_{n}^{(d)}} |\{Z \in M \mid z(Z) \in \zeta^\square\}| \leq c_4 n, \quad n \geq n_0.
\end{align}
For the rest of the proof we fix an $m \in \mathbf{M}$ and an $n_0 \in \N$ such that \eqref{eq:non-trivial-phase-transition:1} and \eqref{eq:non-trivial-phase-transition:2} hold with $m$ in place of $M$.

First we claim that for any $\alpha \in \cA_n(\cG_m)$ there is a $\zeta \in \cA_{\lfloor 3^d(1-c_3)^{-1} n \rfloor}^{(d)}$ such that $\alpha$ is covered by $\zeta^\square$, i.e.\ $Z \subset \zeta^\square$ for each cell $Z \in \alpha$.

Let us assume, one needs a larger lattice animal $\zeta$ to cover the cells of $\alpha$ with $\zeta^{\square}$. Then \eqref{eq:non-trivial-phase-transition:1} implies that there is a subset $\tilde \zeta \subset \zeta$ of more than $3^d n$ vertices $v$ with $U_v = 0$. Hence we would find a subset $\bar \zeta\subset \tilde \zeta$ of more than $n$ vertices, where for any two distinct vertices $v, w \in \bar \zeta$ we have $\|v - w\|_\infty \geq 3$. Each box $v^{\square}$ with $v \in \bar \zeta$ is intersected by a cell of $\alpha$ then, that doesn't intersect any other box of $\bar \zeta$. Thus $\alpha$ would have more than $n$ vertices, which is a contradiction.

If $\alpha$ is covered by $\zeta^\square$, then $|\alpha| \leq \sum_{v \in \zeta} Y_v \leq c_4 \lfloor 3^d(1-c_3)^{-1} n \rfloor$. Hence $|\cA_n(\cG_m)|$ can't be larger than the number of lattice animals of size $\lfloor 3^d(1-c_3)^{-1} n \rfloor$ times the number of subsets of a set with $c_4 \lfloor 3^d(1-c_3)^{-1} n \rfloor$ elements. As $|\cA_n^{(d)}|$ is exponentially bounded \cite{grimmett1999percolation}, this proves the exponential bound \eqref{eq:non-trivial-phase-transition:claim1} and hence $p_c > 0$.

We will go on by proving the bound \eqref{eq:non-trivial-phase-transition:claim2} on the size of the boundary of $\alpha$. Each cell in $\partial^+ \alpha$ intersects $\zeta^\square$. Hence no cell of $\partial^+ \alpha$ intersects a box $v^\square$ with $U_v(m, \delta) = 0$ and $\min_{w \in \zeta}\|w-v\|_\infty > 1$ . If we define $\xi := \{v \in \Z^d \mid \min_{w \in \zeta} \|v - w\|_\infty  \leq 2\}$ the last argument implies, that $\partial^+ \alpha \subset D(\zeta)^\square$ where
\[
  D(\zeta) := \xi \cup \{v \in \Z^d \mid \text{$v \leftrightarrow \xi$ in $U$}\}
\]
and ``$v \leftrightarrow \xi$ in $U$'' means that there is a path $\gamma$ in $\cZ^d$ that connects $v$ with $\xi$ such that $U_w(m,\delta) = 1$ for all $w \in \gamma$. As $D(\zeta)$ is a lattice animal with a size of at least $n$, we obtain
\[
  \sum_{v \in D(\zeta)} U_v \leq c_3 |D(\zeta)|.
\]
Moreover, $U_v = 1$ for each $v \in D(\zeta) \setminus \xi$ and hence
\[
  |D(\zeta)| - |\xi| \leq \sum_{v \in D(\zeta)} U_v.
\]
If follows that
\[
  |D(\zeta)| \leq \frac{|\xi|}{1-c_3} \leq \frac{5^d |\zeta|}{1-c_3}.
\]
By the choice of $\zeta$ and the fact that $\partial^+ \alpha \subset D(\zeta)^\square$ we find the constant $c_2$ such that $|\partial^+ \alpha| \leq c_2 n$.

Finally we adapt Peierls argument to show $p_c < 1$. We define the 2-dimensional sublattice $\Z^2_d := \{(v,0,\dots,0) \in \Z^d \mid v \in \Z^2\}$ and want to show that for large enough $p \in (0,1)$ there is a.s.\ an infinite black cluster in $m_X$ even if we restrict it to the slab $(\Z_d^2)^\square$. Let $W_v := W_v(m)$, $v \in \Z^2_d$ be the event ``$v^\square$ is intersected by a white cell of $m_X$''.

Let us consider the case, where there is a white cut set $S$ in $\cG_{m,X}$, i.e.\ a set of white vertices such that any infinite path starting in the root intersects $S$. If we choose $\tilde S \subset \Z_d^2$ minimal such that $\tilde S^\square$ covers $S$, we obtain a cut set of $\Z_d^2$, i.e.\ a set such that any path in $\cZ^d$ using only vertices of $\Z_d^2$ and starting in the origin intersects $\tilde S$. It is well known, that in the graph $\cZ_d^{2*} := (\Z_d^2, \{\{u,v\} \subset \Z_d^2 \mid \|u - v\|_\infty = 1\})$ any such cut set contains a cycle $\gamma$ around the origin. Hence if there is no infinite black cluster in $\cG_{m,X}$, then for any $n_1 \in \N$ there is a cycle $\gamma$ in $\cZ_d^{2*}$ around the origin with $|\gamma| \geq n_1$ and such that $W_v$ holds for each $v \in \gamma$.

If $|\gamma| = k$ then all elements of the intersection of $\gamma$ and the positive $x$-axis $\{(x,0,0, \dots, 0) \in \Z^d \mid x \in N_0\} \subset \Z^2_d$ have a first coordinate smaller than $k$. Hence there is an animal $\alpha \in \cA_l(\cZ_d^{2*})$ with $k \leq l \leq 2k$ such that $\gamma \subset \alpha$ and we conclude that
\begin{align}
\begin{aligned}\label{eq:non-trivial-phase-transition:3}
  & \quad \ \P_{p}[\text{there is no infinite black cluster in $\cG_{m,X}$}] \\
  & \leq \sum_{n \geq n_1} \sum_{l = n}^{2n} \sum_{\alpha \in \cA_l(\cZ_d^{2*})} \P_{p}\left[ \bigcap_{v \in \alpha} W_v\right].
\end{aligned}
\end{align}

Due to \eqref{eq:non-trivial-phase-transition:1} we have, that for any $\alpha \in \cA_l(\cZ_d^{2*})$, $l \geq n_0$ the set $\beta_1 := \{v \in \alpha \mid U_v = 0\}$ has at least $(1 - c_3) |\alpha|$ elements. Hence we find a set $\beta_2 \subset \beta_1$ such that $|\beta_2| \geq (1 - c_3) |\alpha| / 3^d$ and $\|v - w\|_\infty \geq 3$ for any $v \neq w \in \beta_2$. The choice of $\beta_2$ implies that no cell that intersects $v^\square$ can intersect $w^\square$ if $v \neq w \in \beta_2$. Hence the corresponding events $W_v$ and $W_w$ are independent and
\begin{align}
  \P_{p}\left[ \bigcap_{v \in \alpha} W_v\right] \leq \P_{p}\left[ \bigcap_{v \in \beta_2} W_v \right] = \prod_{v \in \beta_2} \P_{p}[ W_v].
\end{align}
The center of a cell is contained in the cell by definition. This implies, that if $U_v(m, \delta) = 0$, then the number of cells intersecting $v^\square$ is bounded by $\sum_{w: \|w-v\|_\infty \leq 1} Y_w(m,\delta)$ and
\begin{align}\label{eq:non-trivial-phase-transition:4}
  \P_{p}\left[ \bigcap_{v \in \alpha} W_v\right] \leq \prod_{v \in \beta_2} \left( 1 - p^{\sum_{w: \|w-v\|_\infty \leq 1} Y_w} \right).
\end{align}

By the construction of $\beta_2$ we know, that $\bigcup_{v \in \beta_2} \{w \mid \|w-v\|_\infty \leq 1\}$ is an animal in $\cZ_d^{2*}$ with a size of at most $3^d |\alpha|$. It follows from \eqref{eq:non-trivial-phase-transition:2} that
\begin{align}\label{eq:non-trivial-phase-transition:5}
  \sum_{v \in \beta_2} \sum_{w: \|w-v\|_\infty \leq 1} Y_w \leq c_4 3^d |\alpha|.
\end{align}
If we maximize the right hand side of \eqref{eq:non-trivial-phase-transition:4} under the condition \eqref{eq:non-trivial-phase-transition:5}, we obtain a maximum, when all exponents are of equal size. Hence
\begin{align*}
  \P_{p}\left[ \bigcap_{v \in \alpha} W_v\right] \leq \prod_{v \in \beta_2} \left( 1 - p^{\frac{3^{2d} c_4}{1-c_3}} \right) \leq \left( 1 - p^{\frac{3^{2d} c_4}{1-c_3}} \right)^{(1-c_3) |\alpha| 3^{-d}}
\end{align*}
and choosing $p$ large enough gives us that the right hand side of \eqref{eq:non-trivial-phase-transition:3} is summable and hence arbitrary small for large enough $n_1$. This implies, that $p_c(\cG_m)$ is less than some constant smaller than one, that depends only on $c_3, c_4, \delta$ and $d$. An application of \eqref{eq:doubly-stochastic-model} finishes the proof.\qed \bigskip

\subsection*{Examples}
To give examples for tame tessellations we investigate the conditions (T1) and (T2) seperately. Starting with (T1) we observe, that it only depends on the point process of cell centers. The following Lemma will give a sufficient condition for (T1) with the help of the Laplace functional.

\begin{Lemma}\label{le:tame1-sufficient-cond}
  Let $\Phi := \{z(Z) \mid Z \in M\}$ be the point process of cell centres of the random tessellation $M$. If there is a $t > 0$ and $c_1 \in \R$ such that for large enough $n \in \N$
  \begin{align}\label{eq:tame1-suff-cond}
    L_{\Phi}(-t \ind_{\alpha^\square}) = \E\left[ \exp(t \Phi(\alpha^\square)) \right] \leq c_1^n, \quad \alpha \in \cA_n^{(d)},
  \end{align}
  then there is a $c_2 \in \R$ such that (T1) holds, i.e.
  \[
    \limsup_{n \to \infty } \max_{\alpha \in \cA_n^{(d)}} \frac{1}{n} \sum_{v \in \alpha} Y_v \leq c_2.
  \]
\end{Lemma}
\emph{Proof:} We use sub-additivity and the Markov-inequality to obtain that
\begin{align*}
  \P\left[ \max_{\alpha \in \cA_n^{(d)}}\sum_{v \in \alpha} Y_v > c_2 n\right] & \leq \sum_{\alpha \in \cA_n^{(d)}} \P\left[ \sum_{v \in \alpha} Y_v > c_2 n\right] \\
  & \leq \sum_{\alpha \in \cA_n^{(d)}} e^{-c_2 t n} \E\left[ \exp\left( t \sum_{v \in \alpha} Y_v \right) \right] \\
  & \leq \sum_{\alpha \in \cA_n^{(d)}} e^{-c_2 t n} \E\left[ \exp(t \Phi(\alpha^\square)) \right] .
\end{align*}
As $|\cA_n^{(d)}|$ is exponentially bounded \cite[p. 82]{grimmett1999percolation}, the Borel-Cantelli Lemma yields the assertion if $c_2$ is large enough.\qed \bigskip

We remark, that it is easy to verify, that the above lemma does not depend on the specific value of $\delta$ going into $\alpha^\square$.

The first class of point processes, which fulfill the conditions of Lemma \ref{le:tame1-sufficient-cond} are the so called \emph{$\alpha$-weakly sub-Poisson} processes. A point process $\Phi$ is called $\alpha$-weakly sub-Poisson iff
\begin{align}
  \E\left[ \prod_{i = 1}^k \Phi(A_i) \right] \leq \prod_{i = 1}^k \E[\Phi(A_i)]
\end{align}
for all $k \in \N$ and disjoint Borel sets $A_1, \dots, A_k$. Blaczszyczin and Yogeshwaran have shown the following proposition.

\begin{NamedProp}(\cite[Proposition 2.2]{blaszczyszyn2014stochastic})\\
  Let $\Phi$ be a simple point process on $\R^d$ with intensity $\gamma \in (0,\infty )$ and let $f:\R^d \to [0,\infty )$ be measurable. If $\Phi$ is $\alpha$-weakly sub-Poisson, then
  \[
    \E\left[ \exp\left( \int_{\R^d} f(x) \Phi(dx)\right)\right] \leq \exp\left( \gamma \int_{\R^d} (e^{f(x)} - 1) dx\right).
  \]
\end{NamedProp}

Hence if the point process of cell centers is $\alpha$-weakly sub-Poisson the corresponding random tessellation fulfills (T1) as the above proposition implies that $L_\Phi(-t \ind_{\alpha^\square}) \leq e^{\gamma (e^t-1) \delta^d n}$.

In \cite{blaszczyszyn2014stochastic} some examples of $\alpha$-weakly sub-Poisson processes are given. Most prominently is it shown in Example 2.12 that determinantal point processes fall into this class.

Another important point process class that fulfills (T1) are Poisson cluster processes as we defined them before Theorem \ref{thm:poisson-cluster-is-smp}.

\begin{Lemma}
  Let $\Phi$ be a Poisson cluster process. If $\E[e^{t \Psi_1(\R^d)}] < \infty $ for some $t > 0$ then
  \[
    L_\Phi(-t \ind_{\alpha^\square}) \leq \exp\left( \gamma_0 \delta (\E[e^{t \Psi_1(\R^d)}] - 1) n\right).
  \]
  In particular if $\Phi$ is the process of cell centers of a random tessellation $M$ this tessellation fulfills (T1).
\end{Lemma}
\emph{Proof:} The Laplace functional for Poisson cluster processes is well known (see \cite{daley2003introduction}) and hence our starting point is
\[
  L_{\Phi}(-t \ind_{\alpha^\square})  = \exp\left( - \int_{\R^d} \int 1 - e^{t \mu(\alpha^\square-x)} \P_{\Psi_1}(d\mu) \gamma_0 dx \right).
\]
We recall that $1-e^x \geq x(1- e^{y})/y$ for $x \in [0,y]$ and $y \geq 0$ to obtain that
\begin{align*}
  L_{\Phi}(-t \ind_{\alpha^\square}) & \leq \exp\left( - \gamma_0 \int_{\R^d} \int \mu( \alpha^\square - x) \frac{1- e^{t \mu(\R^d) }}{\mu( \R^d) } \P_{\Psi_1}(d\mu) dx \right) \\
  & = \exp\left( - \gamma_0 \int \int_{\R^d} \mu( \alpha^\square - x) dx \frac{1- e^{t \mu(\R^d) }}{\mu( \R^d) } \P_{\Psi_1}(d\mu) \right) \\
  & = \exp\left( - \gamma_0 \delta^d n \int 1- e^{t \mu(\R^d) } \P_{\Psi_1}(d\mu) \right)
\end{align*}
where the third step is due to Fubini.\qed\bigskip

We want to remark, that the bound in the above Lemma is sharp which can be seen by interpreting the usual Poisson process as a Poisson cluster process where each cluster consists of exactly one point at the origin.

A simple calculation shows that Cox processes fulfill \eqref{eq:tame1-suff-cond} if their generating measure satisfies \eqref{eq:tame1-suff-cond}.

Next we turn to condition (T2) which depends much more on the exact construction of the random tessellation rather than only the cell centers. Again we where not able to show that Poisson hyperplane or STIT tessellations fit into this framework. Hence we stick we with Voronoi tessellations again. This gives us the opportunity to state conditions on the generating point process in terms of void probabilities.

\begin{Lemma}\label{le:tame2-suff-cond}
  Let $\Phi$ be a stationary point process on $\R^d$ and let $M := \cV(\Phi)$ be the induced Voronoi tessellation. If there is a $\delta_1 > 0$ and a $c_1 < 1$ such that for all $n \in \N$ and $I \subset \Z^d$ , $|I| = n$ we have
  \begin{align}\label{eq:tame2-suff-cond}
    \P[\Phi(I^{\square\delta_1}) = 0] \leq c_1^n,
  \end{align}
  then $M$ fulfills (T2).
\end{Lemma}
\emph{Proof:} If $J \subset \Z^d$ with $|J| = n$, then there is an $I \subset \Z^d$ with $|I| = 3^d n$ such that $J^{\square 3\delta_1} = I^{\square \delta_1}$ and hence
\begin{align}
  \P[\Phi(J^{\square 3 \delta_1}) = 0] = \P[\Phi(I^{\square \delta_1}) = 0] \leq c_1^{3^d n}.
\end{align}
Hence without loss of generality we may assume, that $c_1$ is arbitrarily small (it will become clear later in the proof how small we need $c_1$ to be).

We choose an odd $c_2 \in \N$ such that $c_2 > 2 \sqrt{d}$ and define $\delta_2 := c_2 \delta_1$. We need a condition that ensures $U_v(M, \delta_2) = 0$ for some $v \in \Z^d$. To this end we cover each box $w^{\square \delta_1}$, $w \in \Z^d$ with $c_2^d$ boxes of sidelength $\delta_1$. If for each $w \in \Z^d$ with $\|w - v\|_\infty = 1$ all small boxes covering $w^{\square \delta_2}$ contain at least one point, the $U_v(M, \delta_2) = 0$. To proof this, let $x,y \in \R^d$ with $\|x - \delta_2 v\| = \tfrac{\delta_2}{2}$ and $\|y - \delta_2 v\| = \tfrac{3 \delta_2}{2}$. Furthermore let $z_x, z_y$ be the centers of the cells that contain $x$ and $y$ respectively. We know that $x$ and $y$ are contained in a small box with diameter $\sqrt{d} \delta_1$. Hence $\|x - z_x\| \leq \sqrt{d} \delta_1$ and the same is true for $y$ and $z_y$. However the distance between $x$ and $y$ is larger than $\delta_2 = 2 \sqrt{d} \delta_1$ and hence $z_x \neq z_y$. This implies that $y$ is not contained in the same cell as $x$ and hence $U_v(M, \delta_2) = 0$ in this case.

Let $\alpha \in \cA_n^{(d)}$ and $\tilde \alpha \in \cA_{3^d c_2^d n}^{(d)}$ such that $\{v \in \Z^d \mid \exists w \in \alpha:\ \|v - w\|_\infty \leq 1\}^{\square \delta_2} \subset \tilde \alpha^{\square \delta_1}$. No empty small box $v^{\square \delta_1}$ with $v \in \tilde \alpha$ can be responsible for more than $3^d$ vertices $w \in \alpha$ with $U_w(M, \delta_2) = 1$. Hence
\begin{align*}
  \P\left[ \sum_{w \in \alpha} U_w(M, \delta_2) > \frac{n}{2}\right] & \leq \P\left[ \sum_{v \in \tilde \alpha} \ind\{\Phi(v^{\square \delta_1}) = 0\} > \frac{n}{2 \cdot 3^d}\right] \\
  & \leq \P\left[ \exists I \subset \tilde \alpha, |I| \geq \frac{n}{2 \cdot 3^d}:\  \Phi(I^{\square \delta}) = 0\right] \\
  & \leq \binom{3^d c_2^d n}{\lceil \frac{n}{2 \cdot 3^d} \rceil} c_1^{\lceil \frac{n}{2 \cdot 3^d} \rceil}.
\end{align*}
It is well known, that $\binom{n}{c_3 n} < c_4^n$ for any $c_3 \in (0,1)$ and suitable $c_4 \in \R$. Hence for small enough $c_1$ we have
\begin{align*}
  \P\left[ \max_{\alpha \in \cA_n^{(d)}}\sum_{v \in \alpha} U_v(M, \delta_2) > \frac{n}{2} \right] \leq \sum_{\alpha \in \cA_n^{(d)}} \P\left[ \sum_{v \in \alpha} U_v(M, \delta_2) > \frac{n}{2}\right] \leq c_5^n
\end{align*}
for some $c_5 < 1$. The Borel-Cantelli lemma finishes the proof.\qed\bigskip

The first important class of point processes that fulfill \eqref{eq:tame2-suff-cond} are so called $\nu$-weakly sub-Poisson processes. These are processes that have smaller void probabilities than a Poisson process with the same intensity. It was mentioned in \cite[example 2.12]{blaszczyszyn2014stochastic} that determinantal point processes have this property. Other examples can be found in \cite{blaszczyszyn2014stochastic} too. It is easy to check \eqref{eq:tame2-suff-cond} for Cox and Poisson cluster processes.

Let $\Phi$ be a Poisson cluster process. Without loss of generality we may assume, that $\Psi_1(\{0\}) = 1$ a.s.\ where $\Psi_1$ was the first cluster (see the proof of Theorem \ref{thm:poisson-cluster-is-smp}). Hence the void probabilities of $\Phi$ are always less or equal to the void probabilities of the underlying Poisson process $\Phi_0$. This implies, that $\Phi$ satisfies \eqref{eq:tame2-suff-cond}.

For Cox processes the condition \eqref{eq:tame2-suff-cond} may be easily translated into a condition on the Laplace functional of the generating random measure.

\section*{Acknowledgements}
This article covers parts of the results of the authors PhD thesis. The author wants to thank G\"unter Last for his support during this time and for the multitude of fruitful discussions.

\bibliography{C:/Arbeit/LaTeX/myBib}
\end{document}